\documentclass[12pt,a4paper]{article}
\usepackage{amsmath}
\usepackage{amsfonts}
\usepackage{amssymb}
\usepackage{fancyhdr}
\newcommand{\datum}{Version 2005.02.12}
\date{May 2003, \datum}

\newcommand{\sesq}[2]{< #1 \, , \, #2>} %
\newcommand{\sobolevorder}[0]{\tilde{s}} %
\newcommand{\h}[1]{H_{#1}} %
\newcommand{\dualh}[1]{H^{'}_{#1}} %
\newcommand{\multh}[1]{\tilde{H}_{#1}} %
\newcommand{\dualmulth}[1]{\tilde{H}^{'}_{#1}} %
\newcommand{\sequiv}[0]{\mathcal{S}}  %
\newcommand{\timeh}[0]{\bar{T}}   %
\newcommand{\wienerp}[2]{W^{#1}_{#2}} %
\newcommand{\wienerq}[2]{\bar{W}^{#1}_{#2}} %
\newcommand{\zcpx}[1]{p_{#1}}  %
\newcommand{\zcpxd}[1]{\bar{p}_{#1}}   %
\newcommand{\prtfpxd}[1]{\bar{V}_{#1}}   %
\newcommand{\vol}[2]{\sigma^{#1}_{#2}} %
\newcommand{\dualvol}[2]{\sigma^{'{#1}}_{#2}} %
\newcommand{\drift}[1]{m_{#1}} %
\newcommand{\ltrans}[1]{\mathcal{L}_{#1}}   %
\newcommand{\prtfs}[0]{\mathsf{P}} %
\newcommand{\sfprtfs}[0]{\mathsf{P}_{sf}} %
\newcommand{\Gammapx}[2]{\Gamma^{#1}_{#2}} %
\newcommand{\fwrate}[1]{f_{#1}} %
\newcommand{\spotrate}[1]{r_{#1}} %
\newcommand{\derprod}[2]{\mathsf{D}^{#1}_{#2}} %
\newcommand{\jm}[1]{\mathfrak{j}^{#1}} %
\newcommand{\processjm}[1]{\mathfrak{J}^{#1}} %
\newcommand{\rvjm}[1]{J^{#1}} %
\newcommand{\mder}[1]{D_{#1}} %
\newcommand{\uequiv}[0]{U}  %
\newcommand{\uequivadj}[0]{U^{*}}  %
\newcommand{\qed}[0]{\textbf{QED}}   %
\title{Bond Market Completeness and Attainable Contingent Claims}
\author{Erik Taflin\footnote{Permanent address: EISTI, Ecole International des Sciences du Traitement de l'Information,
Avenue du Parc, 95011 Cergy, France; taflin@eisti.fr}~\footnote{This article was partially prepared at  CEREMADE, Universit\'e Paris IX - Dauphine,
Place du Mar\'echal-de-Lattre-de-Tassigny, 75775 Paris (Cedex 16), France, taflin@ceremade.dauphine.fr}}
\newtheorem{theorem}{Theorem}[section]
\newtheorem{lemma}[theorem]{Lemma}

\newtheorem{corollary}[theorem]{Corollary}
\newtheorem{remark}[theorem]{Remark}
\newtheorem{example}[theorem]{Example}
\newtheorem{condition}{Condition}
\numberwithin{equation}{section}
\begin{document}
\maketitle
\thispagestyle{empty}
\begin{abstract}
\noindent
A  general class, introduced
in \cite{I.E.-E.T}, of continuous time bond markets driven by a standard cylindrical
Brownian  motion $\wienerq{}{}$ in $\ell^{2},$ is considered.
We prove that there always exist non-hedgeable random variables in the space
$\derprod{}{0}=\cap_{p \geq 1}L^{p}$ and that $\derprod{}{0}$ has a dense subset of attainable elements,
if the volatility operator is non-degenerated a.e. Such results were proved in
\cite{Bj-Ka-Ru97} and \cite{Bj-Ma-Ka-Ru97} in the case of a bond market driven by  finite dimensional B.m.
and marked point processes.
We define certain smaller spaces $\derprod{}{s},$ $s>0$
of European contingent claims, by requiring that the integrand in the martingale representation,
with respect to $\wienerq{}{}$, takes values in weighted $\ell^{2}$ spaces $\ell^{s,2},$
with a power weight of degree $s.$
For all $s > 0,$ the space $\derprod{}{s}$  is  dense in $\derprod{}{0}$ and is
independent of the particular bond price and volatility operator processes.

A simple condition in terms of $\ell^{s,2}$ norms is given on the volatility operator processes,
which implies if satisfied,
that every element in $\derprod{}{s}$ is attainable. 
In this context a related problem of optimal portfolios of zero coupon bonds is solved for general
utility functions and volatility operator processes, provided that the $\ell^{2}$-valued market price of risk process has
certain Malliavin differentiability properties.
\end{abstract}

%
%
%
\section{Introduction}
In this paper we consider the problem of completeness of continuous time markets
of zero-coupon bonds, with arbitrary positive
time to maturity. To fix the ideas, contingent claims will be elements of the space
$\derprod{}{0}=\cap_{p \geq 1}L^{p},$ where the $L^{p}$ spaces are defined with respect
to an apriori given probability measure $P.$
Introducing the zero-coupon markets, we follow the Hilbert space construction %
given in reference \cite{I.E.-E.T}, which permits a unified approach to bond and
stock markets. The zero-coupon price $\zcpx{t},$ at a given time $t,$
is as a function of time to maturity an element of a certain Sobolev space $\h{}$
of continuous functions.
The basic random object in the theory is the price curve
$\zcpx{t} \in \h{}.$ This is an $\infty$-dimensional object and the number of randomsources
influencing its evolution is big
(cf.  end of the introduction of \cite{Bj-Ma-Ka-Ru97} and \cite{Carmona-Tehr}).
Since the dimension of $\h{}$ is countable
infinite, it is natural that the price process $\zcpx{}$ is driven by a countable infinite
number of random sources. Therefore, in \cite{I.E.-E.T}  the evolution of the price in $\h{}$
is given by a diffusion model driven by a countable infinite number of independent standard
Brownian motions (Bm.), i.e. a standard cylindrical Brownian  motion.

From the point of view of sources of randomness, the construction in reference
\cite{I.E.-E.T} is complementary to and generalizes
the  bond market driven by a finite dimensional Bm. and marked point processes,
with possible infinite mark space, introduced in \cite{Bj-Ka-Ru97} and
\cite{Bj-Ma-Ka-Ru97}.
A notable difference is that the price process in \cite{Bj-Ka-Ru97} and \cite{Bj-Ma-Ka-Ru97} takes values
in a Banach space (a certain sup. normed subspace of $C([0,\infty[ \,)$), which in
a general approach requires a more sophisticated stochastic integration theory
than for Hilbert space valued processes.
The completeness of a bond market driven by an infinite number of random sources
was studied in \cite{Bj-Ka-Ru97} and \cite{Bj-Ma-Ka-Ru97} for their jump-diffusion model.
It was proved that such a market with an infinite mark space, is approximately complete
(Th. $6.11$ of \cite{Bj-Ma-Ka-Ru97}),
i.e. the set of hedgeable claims is dense in a certain sense in the set of claims.
However, in general this market is not complete (Prop. $4.7$ of \cite{Bj-Ka-Ru97}).
This was further developed in \cite{DeDonno Pratelli}. The basic reason for this complication,
compared with stock markets, is that the martingale operator (cf. formula $(6.8)$ of \cite{Bj-Ma-Ka-Ru97}),
being the product of the discounted zero-coupon price and the volatility operator,
is a.e. a compact operator defined on an $\infty$-dimensional topological vector space (TVS).
In the case of the Hilbert space model in \cite{I.E.-E.T} the situation is similar.
There the martingale operator
is compact a.e. (see formula (5.3) and Remark 5.1 of \cite{I.E.-E.T}).
The hedging operator, i.e. the adjoint of the martingale operator
is then also compact a.e.
Intuitively, the market can only be complete if the hedging operator is a.e.
surjective which is never the case, for a countable %
infinity of Brownian motions.

The first purpose of this article is to establish rigorously that the bond market
with a usual derivative market
such as $\derprod{}{0}$ cannot be complete, in the case of a countable %
infinity of Brownian motions (Theorem \ref{Th D0 non complete}).
This is in strong contrast with the case of a finite number of random sources,
where this market is complete when the volatility  operator satisfies certain
non-degeneracy conditions (cf. \cite{I.E.-E.T} formula (3.8) and Remark 5.3).

This raises naturally the question of how to generalize the usual concept of a complete market,
tailored for finite dimensional markets, to bond markets. If the martingale operator
has trivial kernel a.e. then $\derprod{}{0}$ has a dense subspace of hedgeable
elements (Theorem \ref{Th approx complete}). However, this does not
give any information on what
the subset of hedgeable elements is.
Roughly this corresponds to an approximately
complete market, introduced in the different context of \cite{Bj-Ka-Ru97} and \cite{Bj-Ma-Ka-Ru97}.
The solution adapted in this article
simply consists of restricting the set of contingent claims to an allowed subspace
$A \subset \derprod{}{0}$ which satisfies:
\begin{equation}
(i) \; A \; \text{\textit{is a locally convex complete TVS}} \;  \text{and} \;
(ii) \; A \; \text{\textit{is dense in}} \; \derprod{}{0}.
\label{Intro 1}%
\end{equation}
Condition $(i)$ permits to study if it is possible to choose the hedging portfolio as a continuous
function of the contingent claim.
Condition $(ii)$  implies that the price (if continuous on $\derprod{}{0}$) of each element in $\derprod{}{0}$
is determined by the price of elements in $A.$
The bond market, is then said to be relatively complete with respect to
the allowed set $A$ of contingent claims or just $A$\textit{-complete}, if all elements in $A$ are attainable.
The idea here is that it should be easy to check whether or not %
a contingent claim  $X$ is in $A.$ If $X \in A$ then $X$ is hedgeable by definition, while if $X \notin A$ then
we can only conclude that there is a sequence, not necessarily bounded, of self-financing portfolios
with terminal value converging to $X.$ Since the portfolio sequence can be unbounded
the approximation scheme is difficult to use in practice and one needs at least a measure of risk,
which permits to pick the best ``approximate portfolio'' in the sequence.

The second purpose of the article is to introduce
spaces $\derprod{}{s},$ $s \geq0$
of  allowed European contingent claims satisfying (\ref{Intro 1}) and sufficiently
large to contain all commonly used derivatives, including those with discontinuous pay-off functions.
The main point in the definition of $\derprod{}{s},$ $s\geq0$ is that the integrand in the
stochastic integral representation of elements in $\derprod{}{s}$
decreases uniformly at a rate given by weighted $\ell^{2}$-spaces with norm
$y \mapsto (\sum_{i \geq 1}(1+i^{2})^{s}(y^{i})^{2})^{1/2}.$
The spaces $\derprod{}{s},$ $s \geq 0$ are independent of the particular bond price
and volatility operator processes and $\derprod{}{s} \subset \derprod{}{s'},$ for  $s' \leq s.$
The third purpose of the article is to give conditions on the volatility operator (Condition \ref{uniform cond sigma})
such that the market is $\derprod{}{s}$-complete for certain $s>0$ (Theorem \ref{th main}).

The forth purpose of the article
is to apply the $\derprod{}{s}$-completeness of the market to the optimal portfolio problem considered
in \cite{I.E.-E.T}. There, the optimal terminal discounted wealth $\hat{X}$ was first found (Th.$3.3$ of \cite{I.E.-E.T})
under general conditions
and then a hedging portfolio $\hat{\theta}$ of $\hat{X}$ was constructed
for certain cases (deterministic volatility Th.$3.8;$ finite number of Bm. Th.$3.6$).
$\hat{X}$ is not always hedgeable, but to cover more general situations where it is hedgeable,
so an optimal portfolio  $\hat{\theta}$ exists, we here impose %
\begin{equation}
(iii) \; A \; \text{\textit{is an algebra under pointwise multiplication}}. %
\label{Intro 2}
\end{equation}
In fact, $\hat{X}$ is a $C^{1}$ function, polynomially bounded together with its derivative,
of $dQ/dP$ for a martingale measure $Q.$  %
Knowing that some claim like $\ln (dQ/dP)$ is attainable, we use the algebraic properties
of $A$ to prove that $\hat{X}$ is also attainable. Now, $\derprod{}{s},$ $s>0$ is not
an algebra (Remark \ref{multiplication Ds}). However, we define a subspace $\derprod{1}{s} \subset \derprod{}{s}$
of once Malliavin differentiable contingent claims, which is seen to be an algebra by
generalizing the use made of the Clark-Ocone representation formula in \cite{I.E.-E.T}.
The $\derprod{1}{s}$-completeness of the market leads to a fairly general solution
of the optimal portfolio problem (\ref{opt prob}) (Theorem \ref{th opt port}).
Reference \cite{Pham 2003} studies the optimal portfolio problem, within (essentially)
the framework of the jump-diffusion model of \cite{Bj-Ma-Ka-Ru97}.
Existence of optimal terminal discounted wealth is established. However, the hedging
problem is only studied in the sense of approximate hedging, so it does not establish
the existence of an  optimal portfolio.

We note (Remark \ref{rmk binary option}) that the spaces $\derprod{}{s}$ are more
appropriate for the study of general hedging problems than $\derprod{1}{s},$ since
the latter do not contain non-Malliavin-differentiable claims, in particular not binary options.
A Malliavin-Clark-Ocone formalism was also adapted recently in reference \cite{Carmona-Tehr},
for the construction of hedging portfolios in a Markovian context, with a Lipschitz
continuous (in the bond price) martingale operator. This guaranties that the
Malliavin derivative of the bond price is proportional to the martingale operator
(formula (30) of \cite{Carmona-Tehr}).
Hedging is then achieved for a restricted class of claims, namely  European claims
being a Lipschitz continuous function in the price of the bond at maturity.

The main results are proved in \S\ref{Proofs} and auxiliary needed results,
difficult to find on suitable form, are proved in Appendix \ref{App}.

\smallskip

\noindent
\textbf{Acknowledgment:} The author would like to thank Ivar Ekeland and Nizar Touzi for fruitful
discussions and the anonymous referees for constructive suggestions.

\section{Zero-coupon markets and portfolios}
\label{ZC markets and portfolios}
Following closely \cite{I.E.-E.T}, we first introduce zero-coupon markets and portfolios.
We consider a continuous time zero-coupon market, with some finite time horizon $\timeh >0$.
At any date $t \in \mathbb{T}=[0,\timeh],$ one can trade zero-coupon bonds with
maturity $t+T,$ where the time \emph{to} maturity $T\in [ 0,\infty [ \;.$

Uncertainty is modeled by a complete filtered probability space
$(\Omega ,P,\mathcal{F},\mathcal{A}),$
where  $\mathcal{A}=\{\mathcal{F}_{t}\;|\;0\leq t\leq \timeh \},$ is a filtration
of the $\sigma $-algebra $\mathcal{F}=\mathcal{F}_{\timeh}.$ The random
sources are given by independent Brownian motions $\wienerp{i}{},$ $i \in \mathbb{N}^{*},$
where  $\mathbb{N}^{*}=\mathbb{N}-\{0\}$. The filtration $\mathcal{A}$ is
generated by the $\wienerp{i}{},$ $i\in \mathbb{N}^{*}.$

We denote by $\zcpx{t}(T)$ the price at time $t$ of a zero-coupon
yielding one unit of account at time $t+T,$ $t \in \mathbb{T},$ $T \geq 0,$ so that $\zcpx{t}(0)=1.$
For a zero-coupon price, which is a strictly positive $C^{1}$ function in the time to maturity, the instantaneous forward
rate contracted at $t \in \mathbb{T}$ for time to maturity $T \geq 0$ is
\begin{equation}
\fwrate{t}(T)=- \frac{1}{\zcpx{t}(T)}\frac{\partial \zcpx{t}(T)}{\partial T},
\label{forward rate}
\end{equation}
in the Musiela parameterization,
the spot interest rate at $t$ is $\spotrate{t}=\fwrate{t}(0)$ and the
discounted zero-coupon price at time $t$ is $\zcpxd{t}=\zcpx{t}\exp (-\int_{0}^{t}\spotrate{\tau} d\tau ).$

We introduce Hilbert spaces $\h{}$ and  $\multh{}$  of continuous real-valued functions,
which will play the role of state spaces
of the price process $\zcpx{}$ and of drift and volatility processes respectively.
Given $\sobolevorder \in  \; ]1/2,1[ \, ,$ let $\h{}$ be the subspace of all $f \in L^{2}([0, \infty [ \,)$
satisfying
\begin{equation}
\int_{x \geq 0}|f(x)|^{2}\,dx
     +\int_{x, \, y \geq 0}|f(x)-f(y)|^{2}|x-y|^{-1-2\sobolevorder}\,dx\,dy < \infty.
\label{A. equiv Hs norm}
\end{equation}
$\h{}$ is a Sobolev space, which we now give its usual Hilbert space structure, often more
easy to use than the one defined by the equivalent norm given by (\ref{A. equiv Hs norm}).
For $s\in \mathbb{R,}$ let $H^{s}$ (cf. \S7.9 of \cite{Horm}) be the usual Sobolev space of real tempered distributions $f$
on $\mathbb{R}$ such that the function $x \mapsto (1+|x|^{2})^{s/2}\hat{f}(x)$ is an element
of $L^{2}(\mathbb{R}),$ where $\hat{f}$ is the Fourier transform\footnote{\label{Fourier}In
           $\mathbb{R}^{n}$ we denote $x \cdot y=\sum_{1 \leq i \leq n} x_{i} y_{i},$ $x,y \in \mathbb{R}^{n}$ and
           we define the Fourier transform $\hat{f}$ of $f$ by
           $\hat{f}(y)=(2 \pi)^{-n/2}\int_{\mathbb{R}^{n}} \exp (-iy \cdot x) f(x)dx.$}
of $f,$ endowed with the norm:
\begin{equation*}
\| f\| _{H^{s}}=( \int (1+|x|^{2})^{s}\;| \hat{f}(x)|^{2}dx)^{1/2}.
\end{equation*}
We note that by  Plancherel's Theorem (cf. \S2, Ch. VI of \cite{Yosida})  $H^{0}=L^{2}(\mathbb{R}).$
The dual $(H^{s})'$ of $H^{s}$ is identified with $H^{-s}$ by the  continuous bilinear form
 $\sesq{}{} :H^{-s}\times H^{s} \mapsto \mathbb{R}$ :
\begin{equation}
\sesq{f}{g}=\int \overline{ \hat{f}(x) }\text{ } \hat{g}(x)dx,
\label{sesqprod}
\end{equation}
where $ \overline{z}$ is the complex conjugate of $z.$
If $s>1/2$ and $f, g \in H^{s},$ then $f$ is H\"older continuous of order $s-1/2,$
$f(x) \rightarrow 0,$ as $|x| \rightarrow \infty$ and there exists a constant
$C$ independent of $f$ and $g$ such that $\| fg\| _{H^{s}} \leq C \| f\| _{H^{s}} \| g\| _{H^{s}}$
(cf. \cite{Calderon} and \S7.9 of \cite{Horm}).
In particular if $s>1/2,$ then $H^{s}\subset C^{0}\cap L^{\infty }.$
We remind that $\sobolevorder>1/2.$
In $H^{\sobolevorder},$  consider the set $H_{-}^{\sobolevorder}$ of functions with support in
$]-\infty,0 ] \,.$  It is a closed subspace of $H^{\sobolevorder}$, so the quotient space
\begin{equation}
\h{}=H^{\sobolevorder}/H_{-}^{\sobolevorder}
\label{H}
\end{equation}
is a Hilbert space as well. It follows  (cf. formula (7.9.4) of \cite{Horm}) that
the norm defined by (\ref{A. equiv Hs norm}) is equivalent to $\|  \, \; \| _{\h{}}.$
To sum up, a real-valued function $f$ on $[ 0,\infty[$ belongs
to $\h{}$ if and only if it is the restriction to $[ 0,\infty[ $ of
some function in $H^{\sobolevorder}$ and that $\dualh{},$ the dual of $\h{},$ is
the set of all distributions in $H^{-\sobolevorder}$ with support
in $[0,\infty [ \,.$ In particular, $\h{}$ inherits from $H^{\sobolevorder}$ the property of
being a Banach algebra and $\dualh{}$  contains all bounded Radon
measures with support in $[0,\infty[ \,.$

To motivate the introduction of $\multh{}$ suppose that $F \in \multh{}$ is the value
of lets say a volatility process. We then impose that $F g \in \h{}$ for all $g \in \h{}.$
As we saw, if $F \in \h{},$ then this condition is satisfied. It is also the case for
functions $F$ on $[0, \infty[ \,,$ such that $F=a+f,$ for some $a \in \mathbb{R}$ and $f \in \h{}.$
Such $F$ permits to consider volatilities, which do not go to zero when the time to maturity
goes to infinity. We here choose $\multh{}$ to be functions of this form, even if
more general spaces are possible. Then the Hilbert space $\multh{}=\mathbb{R} \oplus \h{},$
since the decomposition of $F=a+f,$ $a \in \mathbb{R}$ and $f \in \h{}$ is unique.
The norm is given by
\begin{equation}
\|F \|_{\multh{}}=(a^{2}+\|f\|_{\h{}}^{2})^{1/2}.
\label{norm mutiplier 0}
\end{equation}
The dual $\dualmulth{},$ of $\multh{}$ is identified with $\mathbb{R} \oplus \dualh{}$
by extending the bi-linear form, defined in (\ref{sesqprod}), to $\dualmulth{} \times \multh{}:$
\begin{equation}
\sesq{F}{G}=ab + \sesq{f}{g},
\label{ext-sesqprod}
\end{equation}
where $F=a+f \in \dualmulth{},$ $G=b+g \in \multh{},$ $a,b \in \mathbb{R},$ $f \in \dualh{}$ and $g \in \h{}.$

In order to introduce the bond dynamics, let  $ \ltrans{}$ denote the
semigroup of left translations defined on real functions on $[0,\infty[ \;:$
\begin{equation}
( \ltrans{a}f) ( T) =f( a+T)
\label{lefttr}
\end{equation}
where  $a \geq 0,$  $T \geq 0.$
$ \ltrans{}$ acts as a strongly continuous contraction semi-group in $\h{}$ (resp. $\multh{}$). The infinitesimal generator
is denoted $\partial$ and its domain\footnote{$\mathcal{D}(A),$ $\mathcal{K}(A)$ and $\mathcal{R}(A)$ denote
respectively the domain, the kernel and the range of a linear operator $A.$
$B^{\perp}$ stands for the annihilator of a subset $B$ of a TVS.}
$\mathcal{D}(\partial)$ is denoted $\h{1}$ (resp. $\multh{1}$).
The norm in $\h{1}$ (resp. $\multh{1}$) is defined by
\begin{equation}
\|f \|_{\h{1}}=(\|f\|_{\h{}}^{2}+\|\partial f \|_{\h{}}^{2})^{1/2} \quad
(\text{resp.} \; \|F \|_{\multh{1}}=(\|F \|_{\multh{}}+\|\partial F \|_{\multh{}}^{2})^{1/2}).
\label{norm Hn mhn}
\end{equation}

Throughout the paper, we shall assume that,  $\zcpxd{}$  is a continuous strictly positive
$\h{1}$-valued progressively measurable processes, with respect to
$\mathcal{A},$ given by an equation of the HJM type (see \cite{HJM92} and equation (2.11) of \cite{I.E.-E.T})
\begin{equation}
 \zcpxd{t}= \ltrans{t }\zcpxd{0}+\int_{0}^{t}\ltrans{t-s}(\drift{s}\zcpxd{s})ds
          +\int_{0}^{t}  \sum_{i \in \mathbb{N}^{*}}\ltrans{t-s}(\vol{i}{s}\zcpxd{s}) d\wienerp{i}{s},
 \label{bond dyn integ p*}
\end{equation}
with boundary condition
\begin{equation}
\zcpxd{t}(0)=\exp (\int_{0}^{t} \frac{\partial \zcpxd{s}(0)}{\zcpxd{s}(0)}ds),
 \label{p* bound cond}
\end{equation}
for $t \in \mathbb{T},$ where $\vol{i}{t},$ $i\in \mathbb{N}^{*},$
and $\drift{t}$ are progressively measurable $\multh{}$-valued processes satisfying
\begin{equation}
 \vol{i}{t}(0)=0 \; \text{for} \; i \in \mathbb{N}^{*}
 \label{cond sigma}
\end{equation}
and
\begin{equation}
 \drift{t}(0)=0.
 \label{cond drift}
\end{equation}
Identifying\footnote{\label{l2}$\ell^{2}$ is the usual Hilbert space of all real sequences $x=(x_{1}, \ldots, x_{n}, \ldots),$
with norm $\|x\|_{\ell^{2}}=\sum_{n \geq 1}(x_{n})^{2}.$}
 $\wienerp{}{}$ with
a $\ell^{2}$ cylindrical Wiener process
(c.f.  \S4.3.1 of reference \cite{DaPrato-Zabczyk}), equation (\ref{bond dyn integ p*}) can be written on
a more compact form. Let $\vol{}{}$ be the progressively measurable
$L(\ell^{2},\multh{})$-valued\footnote{\label{lin cont maps} $L(E,F)$ denotes the space
of linear continuous mappings  from $E$ into $F,$ $L(E)=L(E,E).$ } process
defined by $\vol{}{t}x=\sum_{i \in \mathbb{N}^{*}}\vol{i}{t}x_{i}.$ Then equation (\ref{bond dyn integ p*})
reads (cf. equation (5.5) of \cite{I.E.-E.T})
\begin{equation}
 \zcpxd{t}= \ltrans{t }\zcpxd{0}+\int_{0}^{t}\ltrans{t-s}\zcpxd{s}\drift{s}ds
          +\int_{0}^{t}\ltrans{t-s} \zcpxd{s} \vol{}{s} d\wienerp{}{s}.
 \label{bond dyn integ 1 p*}
\end{equation}
We shall assume that $\vol{}{}$ takes its values in the subspace of Hilbert-Schmidt operators of $L(\ell^{2},\multh{}),$
which permits to give a meaning to the stochastic integral in equation (\ref{bond dyn integ 1 p*})
(cf.  \S4.3.1 of \cite{DaPrato-Zabczyk}).

A portfolio is an $\dualh{}$-valued progressively measurable process $\theta$
defined on $\mathbb{T}.$ If $\theta$ is a  portfolio, then its discounted value  at time $t$ is
\begin{equation}
\prtfpxd{t}(\theta)=\sesq{\theta_{t}}{ \zcpxd{t}}.
  \label{wealth *}
\end{equation}
$\theta$ is   an \emph{admissible portfolio} if
\begin{equation}
 \| \theta \|^{2}_{\prtfs}
  =E\left( \int_{0}^{\timeh} (\| \theta_{t} \|^{2}_{\dualh{}}+ \| \dualvol{}{t}\theta_{t}\zcpxd{t} \|^{2}_{\ell^{2}})dt
  +(\int_{0}^{\timeh} | \sesq{\theta_{t}}{\zcpxd{t} \drift{t}} | dt)^{2}
                    \right) < \infty,
  \label{prtf norm}
\end{equation}
where $\dualvol{}{}$ is the adjoint process of $\vol{}{}$ defined  by
$\sesq{f}{\vol{}{t}x}=(\dualvol{}{t}f,x)_{\ell^{2}},$ for all $f \in \dualmulth{}$ and $x \in \ell^{2}.$
Explicitly we have:
\begin{equation}
 \dualvol{}{t}f=(\sesq{f}{\vol{1}{t}},\ldots,\sesq{f}{\vol{i}{t}}, \ldots).
 \label{adjoint p*}
\end{equation}
The set of all admissible portfolios defines  a Banach space $\prtfs$ for the norm $ \| \; \; \|_{\prtfs}.$
A portfolio is self-financing if
\begin{equation}
d\prtfpxd{t}(\theta)=\sesq{\theta_{t}}{\zcpxd{t}\drift{t}}dt
       +\sum_{i \in \mathbb{N}^{*}}\sesq{\theta_{t}}{\zcpxd{t}\vol{i}{t}} d\wienerp{i}{t}.
   \label{SFprtf}
\end{equation}
The  subspace of all self-financing portfolios in $\prtfs$ is  a Banach space $\sfprtfs.$

We next impose a condition on the zero-coupon market.
\begin{condition} \label{market condition}
\text{} \\
a) The initial condition $\zcpxd{0}$ satisfies:
\begin{equation}
\zcpxd{0} \in \h{1} ,\; \zcpxd{0}(0)=1, \; \zcpxd{0} >0,  \label{inital cond}
\end{equation}
b) $\vol{i}{},$ $i\in \mathbb{N}^{*}$  are given progressively measurable $\multh{1}$-valued processes,
such that (\ref{cond sigma}) is satisfied and such that for all $a \in [1,\infty[,$
\begin{equation}
 E( (\int_{0}^{\timeh}\sum_{i \in \mathbb{N}^{*}}\| \vol{i}{t} \|^{2}_{\multh{1}}dt)^{a}
          + \exp(a \int_{0}^{\timeh}\sum_{i \in \mathbb{N}^{*}}\| \vol{i}{t} \|^{2}_{\multh{}}dt)) < \infty,
 \label{domain sigma i}
\end{equation}
c) There exists a family $\{\Gammapx{i}{} \; | \; i\in \mathbb{N}^{*}\}$ of real-valued progressively measurable processes such that
\begin{equation}
\drift{t}+\sum_{i \in \mathbb{N}^{*}} \Gammapx{i}{t} \vol{i}{t}=0,
  \label{rel drift sigma}
\end{equation}
and
 \begin{equation}
 E( \exp(a \int_{0}^{\timeh}\sum_{i \in \mathbb{N}^{*}} | \Gammapx{i}{t} |^{2}dt)) < \infty,
   \quad \forall a\geq 0.
 \label{gamma strong}
\end{equation}
\end{condition}
This condition gives a mathematical meaning to and guarantees the existence of a solution to the mixed initial value and boundary
value problem (\ref{bond dyn integ p*}) and (\ref{p* bound cond})
(see Theorem~$2.1$ and Theorem~$2.2$ of \cite{I.E.-E.T}):
\begin{theorem}\label{Th price}
If Condition \ref{market condition} is satisfied,
then equation (\ref{bond dyn integ p*}) has, in the set of continuous
progressively measurable $\h{}$-valued processes, a unique solution $\zcpxd{}.$
This solution has the following properties: $\zcpxd{}$ is a strictly positive continuous
progressively measurable $\h{1}$-valued processes and the boundary condition (\ref{p* bound cond}) is satisfied.
Moreover, if $u \in [1, \infty[$ and $\bar{q}_{t}=\zcpxd{t}/\ltrans{t}\zcpxd{0}$ then
$\zcpxd{} \in L^{u}(\Omega, P, L^{\infty}(\mathbb{T}, \h{1})),$
$\bar{q}, 1/\bar{q} \in L^{u}(\Omega, P, L^{\infty}(\mathbb{T},\multh{1}))$ and
$\zcpxd{}(0), 1/\zcpxd{}(0) \in L^{u}(\Omega, P, L^{\infty}(\mathbb{T},\mathbb{R})).$
\end{theorem}
Condition \ref{market condition} also guarantees the existence of a martingale measure (see Theorem~$2.8$ and Corollary~$2.10$ of \cite{I.E.-E.T}).
In order to state the result let
\begin{equation}
\xi_{t}=\exp{ \left( -\frac{1}{2}\int_{0}^{t}\sum_{i \in \mathbb{N}^{*}} (\Gammapx{i}{s})^{2}ds
    +\int_{0}^{t}\sum_{i \in \mathbb{N}^{*}}\Gammapx{i}{s} d\wienerp{i}{s}\right)},
 \label{rad-nikod-der}
\end{equation}
where $t \in \mathbb{T}.$ %
\begin{theorem}\label{Th mart meas}
If Condition \ref{market condition} is satisfied, then $\xi$ is a martingale with respect to $(P,\mathcal{A})$
and $\sup_{t \in \mathbb{T}}\xi_{t}^{\alpha} \in L^{1}(\Omega, P)$ for each $\alpha \in \mathbb{R}.$
The measure $Q,$ defined by
$$dQ=\xi_{\timeh}dP,$$
is equivalent to $P$ on $\mathcal{F}_{\timeh}$
and $t \mapsto \wienerq{i}{t}=\wienerp{i}{t}-\int_{0}^{t} \Gammapx{i}{s} ds,$
$t \in \mathbb{T},$ $i \in \mathbb{N}^{*}$ are
independent Wiener process with respect to $(Q,\mathcal{A}).$
If moreover $\theta \in \sfprtfs,$ then $\prtfpxd{}(\theta)$ is a $(Q,\mathcal{A})$-martingale and
$E(\sup_{t \in \mathbb{T}}(\prtfpxd{t}(\theta))^{2}) < \infty$.
\end{theorem}
Thus, when Condition \ref{market condition} is satisfied, the self-financing criteria (\ref{SFprtf}) is equivalent to
\begin{equation}
\prtfpxd{t}(\theta)=\prtfpxd{0}(\theta)+\int_{0}^{t}\sum_{i \in \mathbb{N}^{*}}\sesq{\theta_{s}}{\zcpxd{s}\vol{i}{s}} d\wienerq{i}{s}.
   \label{SFprtf-q}
\end{equation}
The expected value of a random variable $X$ with respect to $Q$ is denoted $E_{Q}(X)$
and $E_{Q}(X)=E(\xi_{\timeh}X).$ \
\section{Contingent claims} \label{cont claim}
In this section we consider contingent claims $X \in L^{p}(\Omega ,P,\mathcal{F}),$
for all $p \geq 1,$ i.e. $X \in \derprod{}{0}=\cap_{ p \geq 1}L^{p}(\Omega ,P,\mathcal{F})$
and we assume that Condition \ref{market condition} is satisfied.
This is a convenient space, since it contains most usually traded contingent claims
and it gives an easy mathematical analysis. It has also a certain invariance with respect to
the probability measure $P,$ which we shall formulate in a slightly more general
context.
For a Banach space $F$ we define the vector space $\derprod{}{0}(F)=\cap_{1 \leq p}L^{p}(\Omega ,P,\mathcal{F},F)$
and we denote $\derprod{}{0}=\derprod{}{0}(\mathbb{R}).$
It is endowed with the topology induced by the countable sequence of seminorms
\begin{equation}
X \mapsto \| X \|_{L^{n}(\Omega ,P,\mathcal{F},F)},
    \label{seminorms}
\end{equation}
$n \in \mathbb{N}^{*}.$ $\derprod{}{0}(F)$ is then a Fr\'echet space.
Replacing $P$ by the martingale measure
$Q$ in the definition of $\derprod{}{0}$ gives the same space, so in the sequel of
this section we shall use $Q:$
\begin{lemma} \label{Lm D0 P Q}
If Condition \ref{market condition} is satisfied and if $F$ is a Banach space, then
$$\derprod{}{0}(F)=\cap_{1 \leq p}L^{p}(\Omega ,Q,\mathcal{F},F)$$
and the topology of $\derprod{}{0}(F)$ is induced by the sequence of seminorms
$X \mapsto \| X \|_{L^{n}(\Omega ,Q,\mathcal{F},F)},$ $n \in \mathbb{N}^{*}.$
\end{lemma}
In general there are non-attainable random variables in the space $\derprod{}{0}$
(see  Theorem \ref{Th D0 non complete}).
In order to obtain complete markets,
we shall therefore restrict the set of allowed contingent claims.
To specify various subspaces of $\derprod{}{0}$ of allowed contingent claims we introduce certain Hilbert spaces
and the isomorphism of square integrable random variables and square integrable progressively measurable processes.

For $s \in \mathbb{R},$ let $\ell^{s,2}$ be the  Hilbert space of real sequences
endowed with the norm
\begin{equation}
\|x\|_{\ell^{s,2}}=(\sum_{i \in \mathbb{N}^{*}}(1+i^{2})^{s}|x^{i}|^{2})^{1/2}.
    \label{ls,2}
\end{equation}
Obviously $\ell^{2}=\ell^{0,2}.$
The operator  $\jm{}$ in $\ell^{s,2},$ with domain $\ell^{s+1,2}$ and given by
\begin{equation}
 (\jm{}x)^{k}=kx^{k}
\label{j}
\end{equation}
is selfadjoint and strictly positive. Obviously, if $a \geq 0$ then the domain of $\jm{a}$ is $\ell^{s+a,2}$
and for all $s_{1},s_{2} \geq 0,$ $\jm{s_{1}+s_{2}}=\jm{s_{1}}\jm{s_{2}}$ (domains included).

Let  $\nu$ be the product of the Lebesgue measure on $\mathbb{T}$ %
and the counting measure on $\mathbb{N}^{*}$ and 
let $L^{2}(\mathbb{T} \times \mathbb{N}^{*})$ be the space of real square integrable functions
with respect to  $\nu.$
We have $L^{2}(\mathbb{T} \times \mathbb{N}^{*}) \cong L^{2}( \mathbb{T},\ell^{2}).$
For $p \geq 2,$ $L^{p}(\Omega,Q, L^{2}( \mathbb{T},\ell^{2}))$ is the Banach space
with norm defined by
\begin{equation}
\|x\|_{L^{p}(\Omega,Q, L^{2}( \mathbb{T},\ell^{2}))}=( E_{Q}(\int_{0}^{\timeh}\|x_{t}\|_{\ell^{2}}^{2}dt)^{p/2})^{1/p}.
    \label{norm Lap}
\end{equation}
$L_{a}^{p}$ denotes the closed subspace of all progressively measurable  elements (modulo equivalence) in
$L^{p}(\Omega,Q, L^{2}( \mathbb{T},\ell^{2})).$
For all $s \geq 0$ and $p \geq 2,$ let $L_{a,s}^{p}=L_{a}^{p} \cap L^{p}(\Omega,Q, L^{2}( \mathbb{T},\ell^{s,2})).$
It is a closed subspace of $L^{p}(\Omega,Q, L^{2}( \mathbb{T},\ell^{s,2}))$ and we give
$L_{a,s}^{p}$ the corresponding Banach space structure.
The operator $\processjm{}$ in $L^{p}(\Omega,Q, L^{2}( \mathbb{T},\ell^{s,2})),$ $p \geq 2$ and $s \in \mathbb{R},$
defined by its domain $L^{p}(\Omega,Q, L^{2}( \mathbb{T},\ell^{s+1,2}))$ and the expression
\begin{equation}
 (\processjm{}x)_{t}(\omega)=\jm{}x_{t}(\omega)
\label{J}
\end{equation}
is a closed operator. For $a \geq 0,$ the fractional power given by
$(\processjm{a}x)_{t}(\omega)=\jm{a}x_{t}(\omega)$ and by its domain $L^{p}(\Omega,Q, L^{2}( \mathbb{T},\ell^{s+a,2})),$
is also a closed operator.  By  restriction $\processjm{}$ defines an operator in $L_{a,s}^{p}$
with domain $L_{a,s+1}^{p},$ which we also denote by $\processjm{}.$ For $p=2,$
$\processjm{}$ is selfadjoint and positive. %

For completeness we state the following result, proved in Appendix \ref{App}: %
\begin{lemma} \label{Lm L2 isomrphism}
The operator $\uequiv$ given by
\begin{equation} \label{isomrphism}
\uequiv(c,x)=c+\sum_{i \in \mathbb{N}^{*}}\int_{0}^{\timeh}x_{t}^{i}d\wienerq{i}{t}.
\end{equation}
is unitary on $\mathbb{R} \oplus L_{a}^{2}$ to $L^{2}(\Omega ,Q,\mathcal{F}).$
Let $\uequivadj : L^{2}(\Omega ,Q,\mathcal{F}) \rightarrow \mathbb{R} \oplus L_{a}^{2}$
be the adjoint of $\uequiv.$ If  $p \geq 2,$
then $\uequiv(\mathbb{R} \oplus L_{a}^{p}) = L^{p}(\Omega ,Q,\mathcal{F}),$
$\uequivadj L^{p}(\Omega ,Q,\mathcal{F}) = \mathbb{R} \oplus L_{a}^{p}$
and the restriction of $\uequiv$ to $\mathbb{R} \oplus L_{a}^{p},$ also denoted by $\uequiv,$ defines a
homeomorphism $\uequiv: \mathbb{R} \oplus L_{a}^{p} \rightarrow L^{p}(\Omega ,Q,\mathcal{F}).$
\end{lemma}
We now transport by unitary equivalence the selfadjoint operator $0 \oplus \processjm{}$ in $\mathbb{R} \oplus L_{a}^{2}$
to a positive selfadjoint operator $\rvjm{}=\uequiv (0 \oplus  \processjm{}) \uequivadj$ in $L^{2}(\Omega ,Q,\mathcal{F}).$
For $X$ in the domain of $\rvjm{}$ we have
\begin{equation} \label{def gen in L2}
\rvjm{}X=\sum_{i \in \mathbb{N}^{*}}\int_{0}^{\timeh}(\processjm{}x)_{t}^{i}d\wienerq{i}{t},
\end{equation}
where $x$ is given by Lemma \ref{Lm L2 isomrphism}.
For $p \geq 2$ and $s \geq 0,$ the vector space $\mathbb{D}_{s}^{p}=\uequiv (\mathbb{R} \oplus L_{a,s}^{p})$ is
according to Lemma \ref{Lm L2 isomrphism} a subspace of $L^{p}(\Omega ,Q,\mathcal{F})$
and $\mathbb{D}_{0}^{p}=L^{p}(\Omega ,Q,\mathcal{F}).$
Due to the corresponding properties of $\processjm{},$  for $a \geq 0,$ it follows that 
$\rvjm{a}$ with domain  $\mathbb{D}_{s+a}^{p}$  is a closed operator in $\mathbb{D}_{s}^{p}.$ 
 $\mathbb{D}_{s}^{p}$  becomes a Banach space for the norm defined by:
\begin{equation}
\|X\|_{\mathbb{D}_{s}^{p}}=(\|X\|_{L^{p}(\Omega,Q)}^{p}+\|\rvjm{s}X\|_{L^{p}(\Omega,Q)}^{p})^{1/p}.
 \label{norm Dsp}
\end{equation}
We then define a decreasing family of Fr\'echet spaces
\begin{equation} \label{Ds}
\derprod{}{s}=\cap_{p \geq 2}\mathbb{D}_{s}^{p},
\end{equation}
$s \geq 0$ of contingent claims. Let $\derprod{}{\infty}=\cap_{s \geq 0}\derprod{}{s}.$
\begin{example} \label{ex binary option}
Let us consider the simple case of a deterministic volatility operator $\vol{}{}$
and a binary option with discounted pay-off $X,$ where $X=0$ for $\zcpxd{\timeh}(T) <K$
and  $X=1$ for $\zcpxd{\timeh}(T) \geq K,$ with $T>0$ and $K>0$ given.
An explicit calculus of the discounted value $E_{Q}(X|\mathcal{F}_{t})$ and It\^o's lemma
give that $X=\uequiv(c,x),$ with $c=E_{Q}(X)$ and
$x_{t}^{i}=g(K/\zcpxd{t}(T+\timeh -t))\vol{i}{t}(T+\timeh -t),$ where $g$
is a continuous bounded function. If
$\int_{0}^{\timeh}\sum_{i \in \mathbb{N}^{*}}i^{2s}\| \vol{i}{t} \|^{2}_{\multh{}}dt < \infty,$
then  $X \in \derprod{}{s}.$
\end{example}
Under certain conditions, $\derprod{}{s}$ will turn out to be  a space of
allowed contingent claims in a complete market if $s$ is sufficiently large.
The space $\derprod{}{s}$ is then satisfactory from the point of view of hedging contingent claims
since it contains commonly (and less commonly) used derivatives, including standard and exotic
bond and interest rate options. We remind that the pay-off for such options can be expressed
as a function of the price (curves) $\zcpxd{}.$
However, if $s>0$ then $\derprod{}{s}$ is not closed under multiplication:
\begin{remark} \label{multiplication Ds}
The space $\derprod{}{s},$ $s > 0,$ is not closed under multiplication.
In fact there exists $X \in \derprod{}{\infty}$ such that for all $s > 0,$  $X^{2} \notin \mathbb{D}_{s}^{2}.$
We now construct such a $X.$ Let $c^{i}=((1+i)^{1/2} \log{(1+i)})^{-1}$ and
$a_{t}=\sum_{i \in \mathbb{N}^{*}}c^{i}\wienerq{i}{t},$ which is well-defined since
$\sum_{i \in \mathbb{N}^{*}}(c^{i})^{2} < \infty.$ We set $X=\int_{0}^{\timeh}a_{t}d\wienerq{1}{t}.$
It follows from (\ref{norm Dsp}), (\ref{def gen in L2}) and the definition of $X$ that
$\|X\|_{\mathbb{D}_{s}^{p}}^{p}=2E_{Q}(|X|^{p}),$ for $s \geq 0$ and $p \geq 2.$
Let $Y_{t}=\int_{0}^{t}a_{s}d\wienerq{1}{s}$ %
and $b=\sup_{t \in \mathbb{T}}|a_{t}|.$
It follows using the Burkholder-Davis-Gundy (BDG) inequalities, that for some constants $C_{p}$ and $C'_{p}:$
$$\|X\|_{\mathbb{D}_{s}^{p}}^{p} \leq C_{p}E_{Q}((\int_{0}^{\timeh}(a_{t})^{2}dt)^{p/2})
\leq C_{p}\timeh^{p/2} E_{Q}(b^{p}) \leq  C'_{p} (\sum_{i \in \mathbb{N}^{*}}(c^{i})^{2})^{p/2} < \infty.$$
One checks that
$$X^{2}=E_{Q}(X^{2})+2 \int_{0}^{\timeh}a_{t}Y_{t}d\wienerq{1}{t}
   +2\int_{0}^{\timeh}(\timeh-t)a_{t}\sum_{i \in \mathbb{N}^{*}}c^{i}d\wienerq{i}{t}.$$
The two first terms on the right hand side are in $\mathbb{D}_{s}^{p}$ for all $s \geq 0$ and $p \geq 2.$
However, if $s > 0$ then $\sum_{i \in \mathbb{N}^{*}}(i^{s}c^{i})^{2}$ diverges,
so the third term on the right hand side is not in $\mathbb{D}_{s}^{2}.$
\end{remark} 
The fact that $\derprod{}{s}$ is not closed under multiplication, is a serious draw back
for the construction of optimal portfolios, such as considered in \cite{I.E.-E.T}. Therefore
we shall introduce a decreasing family of Fr\'echet spaces $\derprod{1}{s},$ $s \geq 0,$
where $\derprod{1}{s} \subset \derprod{}{s}$ and $\derprod{1}{s}$ is an algebra
under multiplication.
The measure $\nu$ is atomless, so the Gaussian process $\{(h,\wienerq{}{}(h)) \, | \, h \in L^{2}(\mathbb{T} \times \mathbb{N}^{*})\},$ where
\begin{equation}
\wienerq{}{}(h)= \sum_{i \in \mathbb{N}^{*}}\int_{0}^{\timeh}h(s,i) d\wienerq{i}{s},
 \label{gaussian process}
\end{equation}
is well-defined (cf. \cite{Nual71}). The Malliavin derivative operator $\mder{},$
is also well-defined on smooth random variables:
\begin{equation}
\mder{(t,i)}X=\sum_{l=1}^{n}f_{l}(\wienerq{}{}(h_{1}), \ldots, \wienerq{}{}(h_{n}))h_{l}(t,i),
\label{m deriv}
\end{equation}
where $X=f(\wienerq{}{}(h_{1}), \ldots, \wienerq{}{}(h_{n})),$ $(t,i) \in \mathbb{T} \times \mathbb{N}^{*},$
$n \in \mathbb{N}^{*},$ $f \in C^{\infty}(\mathbb{R}^{n})$  is polynomially bounded together with all its derivatives
and $f_{l}(x_{1}, \ldots, x_{n})=\frac{\partial f}{\partial x_{l}}(x_{1}, \ldots, x_{n}).$
For all $1\leq p < \infty,$ the linear map in (\ref{m deriv}) defines a closed  linear map, also denoted $\mder{},$
from $L^{p}(\Omega ,Q,\mathcal{F})$ to $L^{p}(\Omega ,Q,\mathcal{F},L^{2}(\mathbb{T} \times \mathbb{N}^{*})),$
with dense domain $\mathbb{D}^{1,p}$ (cf. \cite{Nual71}).
$\mder{t}X$ denotes the $\ell^{2}$ valued random variable defined by the canonical
isomorphism $L^{2}(\mathbb{T} \times \mathbb{N}^{*}) \cong L^{2}( \mathbb{T},\ell^{2}).$
We denote by $\mathcal{D}_{0}$ the subset of
random variables $X$ in (\ref{m deriv}) also satisfying the restrictions that $f$ has compact
support and all the $h_{i}$ are finite sequences.

For $p \geq 2$ and for  $s \geq 0,$ $\mder{}\mathcal{D}_{0} \subset L^{p}(\Omega,Q, L^{2}( \mathbb{T},\ell^{s,2})),$
so  the operator $\processjm{s}\mder{}$ from $L^{p}(\Omega ,Q,\mathcal{F})$
to $L^{p}(\Omega,Q, L^{2}( \mathbb{T},\ell^{2}))$ is densely defined.
It is also closed. In fact, let  $\{X_{n}\}_{n \geq 1} \subset \mathcal{D}_{0}$
converge to $X$ in $L^{p}(\Omega ,Q,\mathcal{F})$ and let $\processjm{s}\mder{}X_{n}=y_{n} \rightarrow y$
in $L^{p}(\Omega,Q, L^{2}( \mathbb{T},\ell^{2})).$ The inverse of $\processjm{s}$
in $L^{p}(\Omega,Q, L^{2}( \mathbb{T},\ell^{2}))$ is a continuous operator.
If $\processjm{s}x_{n}=y_{n}$ and $\processjm{s}x=y,$ then
$x \in L^{p}(\Omega,Q, L^{2}( \mathbb{T},\ell^{s,2})),$ the domain of $\processjm{s},$
and $x_{n} \rightarrow x$ in $L^{p}(\Omega,Q, L^{2}( \mathbb{T},\ell^{2})).$
$\mder{}$ is closed, so $x=\mder{}X$ and therefore $y=\processjm{s}\mder{}X.$
The domain $\mathbb{D}_{s}^{1,p}$ of $\processjm{s}\mder{}$ becomes a real Banach space for the norm defined by
\begin{equation}
\| X \|_{\mathbb{D}_{s}^{1,p}}=(\| X \|_{L^{p}(\Omega,Q)}^{p}
+\|\processjm{s}\mder{} X \|_{L^{p}(\Omega,Q, L^{2}( \mathbb{T},\ell^{2}))}^{p})^{1/p}.
\label{sem norm Dsp}
\end{equation}
We note that $\mathbb{D}_{0}^{1,p}=\mathbb{D}^{1,p}.$
The decreasing family of Fr\'echet spaces of contingent claims $\derprod{1}{s},$
$s \geq 0$ is defined by
\begin{equation} \label{Ds1}
\derprod{1}{s}=\cap_{p \geq 2}\mathbb{D}_{s}^{1,p}.
\end{equation}
One easily checks that multiplication is
continuous in $\derprod{1}{s}$ and that $\derprod{1}{\infty}=\cap_{s \geq 0}\derprod{}{s}$ is dense in $\derprod{1}{0},$
which %
is dense both in  $L^{2}(\Omega ,Q,\mathcal{F})$ and in  $L^{2}(\Omega ,P,\mathcal{F}).$
\begin{remark} \label{rmk binary option}
Due to its algebraic structure the space $\derprod{1}{s}$ %
is suitable for solving optimal portfolio problems. However it
 does not include all
commonly used derivative products. For example the only binary options in $\mathbb{D}^{1,2}$
are the trivial ones $X=1$ and $X=0.$ The space $\derprod{}{s}$ does not have this shortcoming
(see Example \ref{ex binary option}).
This indicates that the spaces $\derprod{}{s}$ and $\mathbb{D}_{s}^{p}$ are more appropriate than the spaces
$\derprod{1}{s}$ and $\mathbb{D}_{s}^{1,p}$ when considering general hedging problems.
\end{remark}
The Clark-Ocone representation (see \cite{Nual71}), used in \cite{I.E.-E.T}
(see \S3.2 and Lemma A.5 of \cite{I.E.-E.T}) in a context of portfolio optimization,
generalizes to %
$\mathbb{D}^{1,2},$ so in particular to  $\derprod{1}{s}:$
\begin{lemma} \label{Lem Clark-Ocone}
If $X \in \mathbb{D}^{1,2}$ and $x$ is given by the isomorphism of Lemma \ref{Lm L2 isomrphism},
then $x_{t}=E_{Q}(\mder{t}X\, | \, \mathcal{F}_{t}),$ $t \in \mathbb{T}.$
\end{lemma}
We omit the proof of this lemma, since it is so similar to the proof of the analog
result for a one dimensional Brownian motion (see Proposition 1.3.5 of \cite{Nual71}).
The space of contingent claims $\mathbb{D}_{s}^{1,p}$ is smaller  than  $\mathbb{D}_{s}^{p}:$
\begin{corollary} \label{cor inclusion}
Let $p \geq 2$ and $s \geq 0.$ Then $\mathbb{D}_{s}^{1,p} \subset \mathbb{D}_{s}^{p},$
$\derprod{1}{s}\subset \derprod{}{s}$ and the inclusion maps are continuous.
Moreover $X \in \derprod{1}{s}$ if and only if $X \in \derprod{}{0}$ and
$\mder{}X \in \derprod{}{0}(L^{2}( \mathbb{T},\ell^{s,2})).$
\end{corollary}
\section{Main results} \label{results}
In order to find a hedging portfolio $\theta \in \sfprtfs$ of a contingent claim $X,$
we have according to (\ref{SFprtf-q}) to solve the equation $\sesq{\theta_{t}}{\zcpxd{t}\vol{i}{t}}=x_{t}^{i},$
for all $i \in \mathbb{N}^{*}$ or equivalently
\begin{equation}
(b_{t}(\omega))'\theta_{t}(\omega)=x_{t}(\omega),
  \label{prtf eq 1}
\end{equation}
a.e. $(t,\omega) \in \mathbb{T} \times \Omega,$
where $b_{t}(\omega)=\zcpxd{t}(\omega)\vol{}{t}(\omega) \in L(\ell^{2}, \h{})$ is the martingale operator at $(t,\omega)$ and $x$ is given by
the martingale decomposition of $X$ in Lemma \ref{Lm L2 isomrphism}.
Let $(b_{t}(\omega))^{*}$ be the adjoint of $b_{t}(\omega)$ with respect to the scalar product in $\h{}.$
Using that the operators $b_{t}(\omega)^{*}$ and $(b_{t}(\omega)^{*}b_{t}(\omega))^{1/2}$
have the same range a.e.  $(t,\omega) \in \mathbb{T} \times \Omega,$ (cf. Lemma \ref{range})
we prove, %
that derivatives in $ L^{p},$ $p \geq 1$ and even derivatives in $\derprod{}{0}$  are
not always attainable:
\begin{theorem} \label{Th D0 non complete}
If condition \ref{market condition} is satisfied, then there exists $X \in \derprod{}{0}$
such that $\prtfpxd{\timeh}(\theta) \neq X$ for all $\theta \in \sfprtfs.$
\end{theorem}
For an example and generalizations, see Remark \ref{bounded non-hedgeable}.
The bond market is approximatively complete in the following sense:
\begin{theorem} \label{Th approx complete}
Let condition \ref{market condition} be satisfied.
$\derprod{}{0}$ has a dense subspace of attainable contingent claims
if and only if the operator $\vol{}{t}(\omega)$ has
a trivial kernel a.e. $(t,\omega) \in \mathbb{T} \times \Omega.$
\end{theorem}
To introduce complete markets, we shall impose a supplementary condition on the volatility
operator. We  now give a motivation of this condition.
The operator  $B_{t}(\omega)=l_{t}\vol{}{t}(\omega),$ where $l_{t}=\ltrans{t}\zcpx{0},$ is a.e. $(t,\omega)$
a Hilbert-Schmidt operator from $\ell^{2}$  to  $\h{},$
when Condition \ref{market condition} is satisfied. Let
\begin{equation}
A_{t}(\omega)=(B_{t}(\omega))^{*}B_{t}(\omega)
  \label{At}
\end{equation}
where $(B_{t}(\omega))^{*}$ is the adjoint of $B_{t}(\omega)$ with respect to the scalar product in $\h{}.$
$A_{t}(\omega)$ is a positive self-adjoint trace-class operator in $\ell^{2}$ a.e. $(t,\omega),$
when Condition \ref{market condition} is satisfied.
In particular the operator $A_{t}(\omega)$
in  $\ell^{2}$ is compact a.e. $(t,\omega),$ so it can not have a bounded inverse.
However it can have an inverse  defined on $\ell^{s,2},$ for some $s>0.$
This simple observation leads us to replace the non-degeneracy condition,
which gives complete markets in the case of a finite number of random sources
(see \cite{I.E.-E.T} formula (3.8) and Remark 5.3), by the following:
\begin{condition} \label{uniform cond sigma}
There exists $s >0$ and $k \in \derprod{}{0},$ such that for all  $x \in \ell^{2}:$
\begin{equation}
 \|x\|_{\ell^{2}} 
 \leq k(\omega)  \|(A_{t}(\omega))^{1/2}x \|_{\ell^{s,2}} \; \text{a.e.} \; (t,\omega) \in \mathbb{T} \times \Omega.
  \label{uniform cond sigma eq}
\end{equation}
\end{condition}
As we will see, if Condition \ref{uniform cond sigma} is satisfied and  $X \in \derprod{}{s},$
then a.e. $(t,\omega)$ the equation
\begin{equation}
  B_{t}^{*}(\omega)\eta_{t}(\omega)=x_{t}(\omega)
  \label{prtf eq 2}
\end{equation}
has a solution in $\h{}$ given by
\begin{equation}
 \eta_{t}(\omega)=S_{t}(\omega)(A_{t}(\omega))^{-1/2}x_{t}(\omega),
  \label{prtf eq 3}
\end{equation}
where $S_{t}(\omega),$ the closure of $B_{t}(\omega)(A_{t}(\omega))^{-1/2},$
is an isometric operator from $\ell^{2}$ to $\h{}.$
Let $\sequiv \in L(\multh{},\dualmulth{})$ be defined by
\begin{equation}
  (f,g)_{\multh{}}=\sesq{\sequiv f}{g},
  \label{S equiv}
\end{equation}
for $f,g \in \multh{}.$
The portfolio $\theta^{1},$ given by
\begin{equation}
 \theta^{1}_{t}=(l_{t}/\zcpxd{t}) \sequiv \eta_{t}
  \label{prtf eq 4}
\end{equation}
then satisfies equation (\ref{prtf eq 1}) and gives the risky part of a
self-financing portfolio $\theta=\theta^{0}+\theta^{1} \in \sfprtfs.$
Here $\theta^{0} \in \prtfs$ is a portfolio of zero-coupon bonds with time to maturity $0$ and
\begin{equation}
 \theta^{0}_{t}=a_{t} \delta_{0},
  \label{prtf eq 5}
\end{equation}
where $\delta_{0}$ is the Dirac measure with support at $0$ and $a$ is the unique real valued process
such that  $\theta$ is self-financing.
Heuristically, this leads to %
the completeness of the market, when the allowed
contingent claims are given by $\derprod{}{s}$ and the conditions \ref{market condition} and \ref{uniform cond sigma}
are satisfied:
\begin{theorem} \label{th main}
If Condition \ref{market condition} and Condition \ref{uniform cond sigma} are satisfied
and if $X \in \derprod{}{s},$ where $s>0$ is given by Condition \ref{uniform cond sigma},
then there exists $\theta \in \sfprtfs$ such that $\prtfpxd{\timeh}(\theta)=X.$
Moreover, one such portfolio is $\theta=\theta^{0}+\theta^{1},$
where $\theta^{0},\theta^{1} \in \prtfs \cap \derprod{}{0}(L^{2}(\mathbb{T},\dualh{}{}))$
are given by formulas (\ref{prtf eq 4}) and (\ref{prtf eq 5}). The linear mappings
$\derprod{}{s} \ni X \mapsto \theta^{i} \in \prtfs \cap \derprod{}{0}(L^{2}(\mathbb{T},\dualh{}{})),$
$i=0,1$ are continuous.
\end{theorem}
This theorem has a converse:
\begin{theorem} \label{th converse main}
Let Condition \ref{market condition} be satisfied and assume that there exist
$s \geq 0$ and $k \in \derprod{}{0},$ such that for all  $x \in \ell^{s,2},$
$ \|(A_{t}(\omega))^{1/2}\jm{s}x \|_{\ell^{2}} \leq k(\omega) \|x\|_{\ell^{2}},$
$\text{a.e.} \; (t,\omega) \in \mathbb{T} \times \Omega.$
Then $\prtfpxd{\timeh}(\theta) \in \derprod{}{s},$  for all
$\theta \in \sfprtfs\cap \derprod{}{0}(L^{2}(\mathbb{T},\dualh{}{})).$
\end{theorem}
We shall apply these results to the optimal bond portfolio problem considered in \cite{I.E.-E.T},
which we now introduce.
The set of all admissible  self-financing portfolios with initial wealth $x$ is
\begin{equation*}
\mathcal{C}( x) =\{ \theta \in \sfprtfs \mid \prtfpxd{0}(\theta)=x\}.
\end{equation*}
The optimization problem is, for a given initial wealth $K_{0} ,$ to find a solution
$\hat{\theta} \in \mathcal{C}(K_{0})$ of
\begin{equation}
E(u( \prtfpxd{\timeh}( \hat{\theta})))
   =\sup_{\theta \in \mathcal{C}(K_{0})}E(u( \prtfpxd{\timeh}( \theta)))
  \label{opt prob}
\end{equation}
where the utility function $u$ satisfies the following Inada-type  condition:
\begin{condition} \label{u cond}
\text{} \\
\noindent a) $u:\mathbb{R} \rightarrow \mathbb{R} \cup \{ -\infty \}$ is strictly concave,
upper semi-continuous and finite on an interval $]\underline{x},\infty [,$
with $\underline{x}\leq 0$ (the value $\underline{x}=-\infty $ is allowed). \\ \vspace{-2ex}

\noindent b) $u$ is $C^{2}$ on $]\underline{x},\infty [$
and $u^{\prime }( x) \rightarrow \infty$ when $x\rightarrow \underline{x}$
in $]\underline{x},\infty [.$ \\ \vspace{-2ex}

\noindent c) there exist some $q>0$ and $C>0$ such that %
\begin{equation}
 \liminf_{x \downarrow  \underline{x}} (1+|x|)^{-q}u'(x)>0
  \label{utility cond 1}
\end{equation}
and such that, if $u'>0$ on $]\underline{x},\infty [$ then
\begin{equation}
 \limsup_{x \rightarrow \infty} x^{q}u'(x)< \infty
   \quad \text{and} \quad | x\varphi'(x) | \leq C (x^{q}+x^{-q}) \; \text{for all} \; x>0
  \label{utility cond 2}
\end{equation}
and if $u'$ takes the value zero  then
\begin{equation}
 \limsup_{x \rightarrow \infty} x^{-q}u'(x)<0
   \quad \text{and} \quad  | x\varphi'(x)| \leq C (1+|x|)^{q} \; \text{for all} \; x \in \mathbb{R},
  \label{utility cond 3}
\end{equation}
where $\varphi$ is the inverse of $u'$ restricted to $]\underline{x},\infty [.$
\end{condition}
\begin{theorem}\label{th opt port}
Let Condition \ref{market condition}, Condition \ref{uniform cond sigma} and Condition \ref{u cond} be satisfied
and let $\ln(\xi_{\timeh}) \in \derprod{1}{s},$ where $s>0$ is given by Condition \ref{uniform cond sigma}.
If $K_{0} \in \; ]\underline{x},\infty [\,,$ then problem (\ref{opt prob}) has a solution $\hat{\theta}.$
\end{theorem}
We end this section with the following remarks:
\begin{remark} \label{bounded non-hedgeable}
One can consider stronger formulations of Theorem \ref{Th D0 non complete}. For example,
whether or not one can choose the non-attainable claim $X$ bounded and smooth is an open
question in general. For constant deterministic volatility operators the answer is yes.
In fact, the operator $A_{t}(\omega)=A$ (see (\ref{At})) in $\ell^{2}$
is then constant $(t, \omega)$ and since $A^{1/2}$ is compact, we can choose
$e \in \ell^{2},$ such that $\|e\|_{\ell^{2}}=1$ and $e \notin \mathcal{R}(A^{1/2}).$
Let $g \in C^{\infty}(\mathbb{R})$ be rapidly decreasing together with all its derivatives,
$g(x) > 0$ for all $x \in \mathbb{R},$  $f(x)=\int_{y\leq x}g(y)dy,$ $Y=\sum_{n \geq 1}e^{n} \wienerq{n}{\timeh}$
and $X=f(Y).$ Then, by Lemma \ref{Lem Clark-Ocone},
$X=E_{Q}(X)+\int_{0}^{\timeh}\sum_{n \geq1}x^{n}_{t}d\wienerq{n}{t},$
where $x_{t}=E_{Q}(g(Y)\, | \, \mathcal{F}_{t})e.$ Since  $e \notin \mathcal{R}(A^{1/2})$
and $E_{Q}(g(Y)\, | \, \mathcal{F}_{t}) >0$ for all $t,$ it follows that $X$ is not attainable.
\end{remark}
\begin{remark} \label{weaker conditions}
 Our choice $\ell^{s,2},$ $s \geq 0$  of weighted $\ell^{2}$-spaces, leading to the results of this
section, can be generalized to other weighted $\ell^{2}$-spaces. If $q_{s}$ are the
corresponding norms, then the crucial property
which shall be satisfied is Condition \ref{uniform cond sigma}
with the $\ell^{s,2}$-norms replaced by $q_{s},$  which can depend on $(t,\omega).$
\end{remark}
\begin{remark}\label{rmk opt port}
Conditions can be given directly on $\Gammapx{}{},$ which guarantees that $\ln(\xi_{\timeh})$
satisfies the hypothesis of Theorem \ref{th opt port}. One possibility is:
If $s\geq0$ and if for all $n \in \{0,1,2\}$ and $p \geq 2,$
$\mder{}^{n}\Gammapx{}{} \in L^{p}(\Omega,Q, \otimes^{n+1} L^{2}( \mathbb{T},\ell^{s,2})),$ then
$\ln(\xi_{\timeh}) \in \derprod{1}{s}.$ %
\end{remark}
\section{Proofs}
\label{Proofs}
\noindent \textbf{Proof of Lemma \ref{Lm D0 P Q}}
Let $p \geq 2.$ It follows from Schwarz inequality that
$\|X\|_{L^{p}(P,F)}=(E_{Q}(\xi_{\timeh}^{-1}\|X\|_{F}^{p}))^{1/p} \leq (E_{Q}(\xi_{\timeh}^{-2}))^{1/2p} \|X\|_{L^{2p}(Q,F)}.$
Similarly, $\|X\|_{L^{p}(Q,F)}$ $ \leq (E(\xi_{\timeh}^{2}))^{1/2p} \|X\|_{L^{2p}(P,F)}.$
According to Theorem \ref{Th mart meas}, $E_{Q}(\xi_{\timeh}^{-2})=E(\xi_{\timeh}^{-1}) < \infty$
and $E(\xi_{\timeh}^{2}) < \infty.$
\qed \\ \vspace{-2ex}

\noindent \textbf{Proof of Corollary \ref{cor inclusion}}
Let $p$ and $s$ be as in the corollary and
let $X \in \mathbb{D}_{s}^{1,p}.$ Obviously $X \in \mathbb{D}^{1,2},$ so
we can apply Lemma \ref{Lem Clark-Ocone} giving  $X=\uequiv(c,x),$ where
$x_{t}=E_{Q}(\mder{t}X\, | \, \mathcal{F}_{t}),$ $t \in \mathbb{T}.$ We obtain
$\|\processjm{s}x \|_{L^{p}(\Omega,Q, L^{2}( \mathbb{T},\ell^{2}))}$
$\leq \|\processjm{s} \mder{t}X \|_{L^{p}(\Omega,Q, L^{2}( \mathbb{T},\ell^{2}))}$
$\leq \|X \|_{\mathbb{D}_{s}^{1,p}}.$ Then $x \in L_{a,s}^{p},$ so by definition $X \in \mathbb{D}_{s}^{p},$
proving tat $\mathbb{D}_{s}^{1,p} \subset \mathbb{D}_{s}^{p}.$
The inclusion map is continuous since, by the last inequality and Lemma \ref{Lm L2 isomrphism}, for some
constant $C_{p}$ we have
$\|\rvjm{s}X\|_{L^{p}(\Omega,Q)} \leq C_{p} \|\processjm{s}x \|_{L^{p}(\Omega,Q, L^{2}( \mathbb{T},\ell^{2}))}
\leq C_{p} \|X \|_{\mathbb{D}_{s}^{1,p}}.$

The continuous inclusion of $\derprod{1}{s}$
into  $\derprod{}{s}$ is now a direct consequence of the continuous inclusion of
$\mathbb{D}_{s}^{1,p}$ into $\mathbb{D}_{s}^{p},$ for all $p \geq 2.$

According to formula (\ref{sem norm Dsp}), $X \in \derprod{1}{s}$
if and only if $ X \in L^{p}(\Omega,Q)$ and
$\processjm{s}X \in L^{p}(\Omega,Q, L^{2}( \mathbb{T},\ell^{2})),$ for all $p \geq 2.$
Then by Lemma \ref{Lm D0 P Q}, $X \in \derprod{1}{s}$
if and only if $X \in \derprod{}{0}$ and
$\mder{}X \in \derprod{}{0}(L^{2}( \mathbb{T},\ell^{s,2})).$  \qed  \\ \vspace{-2ex}

\noindent \textbf{Proof of Theorem \ref{Th D0 non complete}}
Using the definition (\ref{S equiv}) of $\sequiv,$ one readily verifies that the
operator $(b_{t}(\omega))'$ in equation (\ref{prtf eq 1}) satisfies
\begin{equation}
(b_{t}(\omega))'\sequiv=(b_{t}(\omega))^{*},
  \label{adj rel}
\end{equation}
a.e. $(t,\omega).$
Since $\sequiv$ defines a homeomorphism of $\dualh{}$ onto $\h{},$ equation (\ref{prtf eq 1})
is equivalent to find a $\h{}$-valued process $y$ satisfying
$\sequiv y \in \sfprtfs$ and $(b_{t}(\omega))^{*}y_{t}(\omega)=x_{t}(\omega),$
a.e. $(t,\omega) \in \mathbb{T} \times \Omega.$
A necessary condition for the existence of a solution of equation (\ref{prtf eq 1})
is then according to Lemma \ref{range}:
\begin{equation}
         x_{t}(\omega) \in \mathcal{R}(((b_{t}(\omega))^{*}b_{t}(\omega))^{1/2}),
       \; \text{a.e.} \; (t,\omega) \in \mathbb{T} \times \Omega.
  \label{prtf eq 6}
\end{equation}
Let $K_{t}(\omega)=((b_{t}(\omega))^{*}b_{t}(\omega))^{1/2},$ let $\|K_{t}(\omega)\|_{H-S}$
be its Hilbert-Schmidt norm and let $\{u_{n}\}_{n \geq 1}$ be the standard orthonormal basis in $\ell^{2}.$ Then
$\|K_{t}(\omega)\|_{H-S}^{2}=\sum_{i \geq 1}\|b_{t}(\omega)u_{i}\|_{\h{}}^{2}
 \leq C \sup_{s} \|\zcpxd{s}(\omega)\|_{\h{}}^{2}  \sum_{i}\|\vol{i}{t}(\omega)\|_{\multh{}}^{2},$
which is integrable according to H\"older's inequality, inequality (\ref{domain sigma i})
of Condition \ref{market condition} and Theorem \ref{Th price}.
Therefore $\|K_{t}(\omega)\|_{H-S}^{2}$ is finite a.e. $(t,\omega).$ We can now apply
Lemma \ref{g_{0} g_{1}} and choose $x_{t}(\omega)=g_{1}(K_{t}(\omega)),$
a.e. $(t,\omega) \in \mathbb{T} \times \Omega.$ Then $x_{t}(\omega) \notin \mathcal{R}(K_{t}(\omega))$
a.e. $(t,\omega) \in \mathbb{T} \times \Omega.$
Since $g_{1}$ is Borel measurable
 on the space of selfadjoint Hilbert-Schmidt operators, with the operator norm topology,
it follows that $x$ is progressive.

Let $X=U(0,x),$ where $U$ is as in Lemma \ref{Lm L2 isomrphism}. It follows from Lemma \ref{g_{0} g_{1}},
Lemma \ref{Lm D0 P Q} and Lemma \ref{Lm L2 isomrphism} that $x \in \derprod{}{0}.$
Since condition (\ref{prtf eq 6}) is not satisfied, it follows that equation (\ref{prtf eq 1})
does not have a solution for this $x.$
\qed \\ \vspace{-2ex}

In the proof of the next two theorems we shall use the following
\begin{lemma} \label{Lm sf prtf}
If Condition \ref{market condition} is satisfied and if
$\theta^{1} \in \derprod{}{0}(L^{2}(\mathbb{T},\dualh{}{}))$ and
$x \in \derprod{}{0}(L^{2}(\mathbb{T},\ell^{s,2}))$ are progressively measurable processes satisfying
formula (\ref{prtf eq 1}),
then $\theta^{1} \in \prtfs.$ If, moreover
$c \in \mathbb{R},$
\begin{equation} \label{proof th main 2}
Y_{t}=c+\sum_{i \in \mathbb{N}^{*}}\int_{0}^{t}x_{s}^{i}d\wienerq{i}{s},
\end{equation}
\begin{equation}
 a_{t}= (\zcpxd{t}(0))^{-1} (Y_{t}-\sesq{\theta^{1}_{t}}{\zcpxd{t}})\; \text{and} \; \theta^{0}_{t}=a_{t} \delta_{0},
 \label{proof th main 5}
\end{equation}
for $t \in \mathbb{T}$ and  $Z=\sup_{t \in \mathbb{T}}|Y_{t}|,$ then
$Y$ is a $Q$-martingale, $Z \in \derprod{}{0},$ $a \in \derprod{}{0}(L^{2}(\mathbb{T})),$
$\theta^{0} \in \prtfs \cap \derprod{}{0}(L^{2}(\mathbb{T},\dualh{}{}))$ and
$\theta=\theta^{0}+\theta^{1} \in \sfprtfs \cap \derprod{}{0}(L^{2}(\mathbb{T},\dualh{}{})).$
The linear map $(\theta^{1},c,x) \mapsto a \in \derprod{}{0}(L^{2}(\mathbb{T}))$ is continuous on
the subspace of progressively measurable processes in
$\derprod{}{0}(L^{2}(\mathbb{T},\dualh{}{})) \times \mathbb{R} \times \derprod{}{0}(L^{2}(\mathbb{T},\ell^{s,2}))$
satisfying (\ref{prtf eq 1}).
\end{lemma}
\textbf{Proof:} Since $\theta^{1}$ satisfies $b'\theta^{1} \in \derprod{}{0}(L^{2}(\mathbb{T},\ell^{s,2}))$
by construction, it follows from the definition of $\prtfs,$ Condition \ref{market condition},
$\theta^{1} \in \derprod{}{0}(L^{2}(\mathbb{T},\dualh{}{})),$
and H\"older's inequality that $\theta^{1} \in \prtfs.$

Obviously $Y$ is a $Q$-martingale. By Doob's $L^{p}$-inequality (cf. \cite{Revuz-Yor}) we have
$\|Z\|_{L^{p}(Q)} \leq c_{p}\sup_{t \in \mathbb{T}}\|Y_{t}\|_{L^{p}(Q)},$
for $p \geq 2.$ We then obtain
$\|Z\|_{L^{p}(Q)} \leq c_{p}\|X\|_{L^{p}(Q)},$
since $|Y|^{p}$ is a $Q$-submartingale. Lemma \ref{Lm D0 P Q} now gives that
\begin{equation} \label{proof th main 3}
Z \in \derprod{}{0}.
\end{equation}
Schwarz inequality and the definition of $Z$ give
$$(\int_{0}^{\timeh}|a_{t}|^{2}dt)^{1/2} \leq (\sup_{t' \in \mathbb{T}}|\zcpxd{t'}(0)|^{-1})
(Z \timeh+ (\sup_{t' \in \mathbb{T}}\|\zcpxd{t'} \|_{\h{}{}}) (\int_{0}^{\timeh}\|\theta^{1}_{t}\|_{\dualh{}{}}^{2}dt)^{1/2} ).$$
Formula (\ref{proof th main 3}), $\theta^{1} \in \derprod{}{0}(L^{2}(\mathbb{T},\dualh{}{}))$ %
and  H\"older's inequality then give
\begin{equation} \label{proof th main 6}
a \in \derprod{}{0}(L^{2}(\mathbb{T}))
\end{equation}
and the announced continuity property of $a.$
Since $\|\delta_{t}\|_{\dualh{}{}} =C < \infty,$ for $t \in \mathbb{T}$ and $C$ independent
of $t,$ it follows from formulas (\ref{cond sigma}) and (\ref{rel drift sigma}) that
$\|\theta^{1}\|_{\prtfs}=C\|a\|_{L^{2}(\Omega ,P,L^{2}(\mathbb{T}))}.$ Formula (\ref{proof th main 6})
now shows that $\theta^{0} \in  \prtfs \cap \derprod{}{0}(L^{2}(\mathbb{T},\dualh{}{})).$

By the definition (\ref{proof th main 5}) of $a,$ it follows that $\prtfpxd{t}(\theta)=Y_{t}.$ $\theta$ is then
self-financing according to formulas (\ref{SFprtf-q}) and (\ref{proof th main 2}) with
$\prtfpxd{0}(\theta)=c.$
\qed \\ \vspace{-2ex}

The following notations will be used in the proof of %
Theorem \ref{Th approx complete}: $d\mu=dtdQ.$
$F$ is the closed subspace of progressively measurable processes in $L^{2}(\Omega,Q, L^{2}( \mathbb{T},\h{})).$
For $p >1,$ the operator $b_{(p)}$ from $L_{a}^{p}$ to $F$ is defined by its domain
$\mathcal{D}(b_{(p)})=
\{x \in L_{a}^{p} \, | \, \int_{\mathbb{T} \times \Omega }\|b_{t}(\omega)x_{t}(\omega)\|_{\h{}}^{2} d\mu < \infty\}$
and
\begin{equation} \label{op b_p}
(b_{(p)}x)_{t}(\omega)=b_{t}(\omega)x_{t}(\omega), %
\end{equation}
where $b$ is as in (\ref{prtf eq 1}).
The operator $b_{(p)}$ is densely defined and closed. %
We note that $b_{(p)},$ $p>1$ is a maximal operator in the sense that it does not
have a nontrivial extension, satisfying (\ref{op b_p}).
The adjoint $(b_{(p)})^{*}$ is given by $((b_{(p)})^{*}y)_{t}(\omega)=(b_{t}(\omega))^{*}y_{t}(\omega),$
$1/p+1/q=1$ and $\mathcal{D}((b_{(p)})^{*})=
\{y  \in F \, | \, E_{Q}((\int_{\mathbb{T}}\|(b_{t}(\omega))^{*}y_{t}(\omega)\|_{\h{}}^{2})^{q/2}) < \infty\}.$
Given a selfadjoint operator $A,$  we denote by $e_{A}$ be the resolution of the identity associated with  $A.$  \\ \vspace{-2ex}

\noindent \textbf{Proof of Theorem \ref{Th approx complete}}
Let  $\mathcal{U}=\{(t,\omega) \in \mathbb{T} \times \Omega \, | \, \mathcal{K}(\vol{}{t}(\omega))\neq \{  0\} \}.$

$1)$ Let $\mu(\mathcal{U})=0.$ The set $\mathcal{D}_{1}= \cap_{p>1} \mathcal{D}^{(p)},$ where
 $\mathcal{D}^{(p)}=\{y \in F \, | \,
    y \in \derprod{}{0}(L^{2}(\mathbb{T},\h{})), b_{(p)}^{*}y \in \derprod{}{0}(L^{2}(\mathbb{T},\ell^{2}))\}$
is dense in $F.$ In fact its enough to consider progressive $y \in L^{\infty}(\Omega \times \mathbb{T},\h{})).$

For $y \in \mathcal{D}_{1},$ we define $\theta^{1}=\sequiv y.$  The relation (\ref{adj rel}) gives
$(b_{t}(\omega))'\theta^{1}_{t}(\omega)=(b_{(p)}^{*}y)_{t}(\omega).$ Then, according to the definition
of $\mathcal{D}_{1},$ the hypotheses of Lemma \ref{Lm sf prtf} are satisfied with
$x_{t}(\omega)=(b_{t}(\omega))'\theta^{1}_{t}(\omega).$ $\theta^{0}$ is defined by (\ref{proof th main 5}).
Lemma \ref{Lm sf prtf} gives $\theta^{0},\theta^{1} \in \prtfs$ and $\theta=\theta^{0}+\theta^{1} \in \sfprtfs.$
Let $\mathcal{D}_{2} \subset \sfprtfs$ be the set of all such $\theta,$ for $y \in \mathcal{D}_{1}.$
By construction $\bar{V}_{\timeh}(\theta)= \uequiv (\bar{V}_{0}(\theta), x)$ so
$\theta \mapsto (\bar{V}_{0}(\theta), x)$ defines a mapping of  $\mathcal{D}_{2}$
onto $\mathbb{R} \oplus  b_{(p)}^{*}\mathcal{D}_{1}.$ $\bar{V}_{\timeh}(\theta) \in \derprod{}{0}$
according to Lemma \ref{Lm L2 isomrphism}, since $x \in \derprod{}{0}(L^{2}(\mathbb{T},\ell^{s,2})).$

For the moment suppose that, for every $p>1,$ $(b_{(p)})^{*}\mathcal{D}_{1}$ is dense in $L_{a}^{q},$ the dual of
$L_{a}^{p},$  where $p^{-1}+q^{-1}=1.$ Since the set $\mathcal{D}_{3}=(b_{(p)})^{*}\mathcal{D}_{1}$ %
is independent of $p,$ it follows that $\mathcal{D}_{3}$ is dense in $\derprod{}{0}(L^{2}(\mathbb{T},\ell^{s,2})).$
By Lemma \ref{Lm L2 isomrphism} it then follows that the set of attainable claims $\uequiv(\mathbb{R} \oplus  \mathcal{D}_{3})$
is dense in $\derprod{}{0}.$

It remains to prove that $(b_{(p)})^{*}\mathcal{D}_{1}$ is dense in $L_{a}^{q}.$
Let $c_{(p)}$ be the restriction of $(b_{(p)})^{*}$ to $\mathcal{D}_{1}$ and let
$\tilde{b}_{(p)}=(c_{(p)})^{*}.$
Then $\tilde{b}_{(p)}$ is an extension of $b_{(p)}=(b_{(p)})^{**}$ and it satisfy (\ref{op b_p})
with $\tilde{b}_{(p)}$ instead of $b_{(p)},$ which shows that $b_{(p)}=\tilde{b}_{(p)}.$
Therefore $\mathcal{K}(b)=\mathcal{R}(c)^{\perp}.$
Since $\mathcal{K}(b)$ is trivial, $\mathcal{R}(c)$ is dense in $L_{a}^{q}.$

$2)$  Let $\mu(\mathcal{U})>0.$ We proceed as in the proof of Theorem \ref{Th D0 non complete}, introduce
the selfadjoint operator
$K_{t}(\omega)=((b_{t}(\omega))^{*}b_{t}(\omega))^{1/2}$ in $\ell^{2},$ apply
Lemma \ref{g_{0} g_{1}} and choose $x_{t}(\omega)=g_{0}(K_{t}(\omega)).$ Then $0 \neq x \in L_{a}^{p}$ for all $p>1$
and $X=U(0,x) \in \derprod{}{0}.$
 $K$ with domain $\mathcal{D}(b_{(2)})$ is selfadjoint in $L_{a}^{2}.$
One readily verifies that $x$ is orthogonal to $\mathcal{R}(K)$ in $L_{a}^{2}.$
This proves that $X=U(0,x)$ is orthogonal to the image of $\mathcal{R}(K)$ under $U(0,\cdot),$
so $L_{a}^{2}$ does not have a  dense subspace of attainable elements. This is then also
the case of $\derprod{}{0}.$
\qed \\ \vspace{-2ex}

\noindent \textbf{Proof of Theorem \ref{th main}}
Let the conditions of the theorem be satisfied. Then $X=\uequiv(c,x)$ for some
$c \in \mathbb{R}$ and $x\in \cap_{p \geq 2}L_{a,s}^{p},$ according to the construction of $\derprod{}{s}.$
We choose  $x $ progressively measurable by changing it on a set of zero measure.
We observe that
$\cap_{p \geq 2}L_{a,s}^{p} \subset \cap_{p \geq 2}L^{p}(\Omega ,Q,L^{2}(\mathbb{T},\ell^{s,2}))
           \subset \derprod{}{0}(L^{2}(\mathbb{T},\ell^{s,2})),$ where the last
relation follows from  Lemma \ref{Lm D0 P Q}. This shows that
\begin{equation} \label{proof th main 1}
x \in \derprod{}{0}(L^{2}(\mathbb{T},\ell^{s,2})),
\end{equation}
where $x$ is progressively measurable.

Let $(t,\omega) \in \mathbb{T} \times \Omega$ be such that $B_{t}(\omega) \in L(\ell^{2},\h{}).$ %
Inequality (\ref{uniform cond sigma eq}) implies
that $(A_{t}(\omega))^{1/2} \in L(\ell^{2})$ has a trivial kernel. Lemma \ref{range}
then firstly shows that $B_{t}(\omega)$ also has a trivial kernel and secondly shows
that $(A_{t}(\omega))^{-1/2}$ is densely defined and that $S_{t}(\omega)$ in formula
(\ref{prtf eq 3}) is isometric from $\ell^{2}$ to $\h{}.$
According to inequality (\ref{uniform cond sigma eq}),
if $z \in \ell^{s,2}$ then $z$ is in the domain of $(A_{t}(\omega))^{-1/2}$ and
$ \|(A_{t}(\omega))^{-1/2}z\|_{\ell^{2}}
 \leq k_{t}(\omega)  \|z \|_{\ell^{s,2}}.$
By equation (\ref{prtf eq 3}) we get
$ \|\eta_{t}(\omega)\|_{\ell^{2}}
 \leq k_{t}(\omega) \|x_{t}(\omega) \|_{\ell^{s,2}}.$ Since this is true a.e. $(t,\omega)$
it follows by integration with respect to $P,$ from Condition \ref{uniform cond sigma} and  H\"older's
inequality, that $\eta \in \derprod{}{0}(L^{2}(\mathbb{T},\h{})).$ $\eta$ is progressively measurable
since this is the case of $x$ and $\vol{}{}.$
In fact, if $y$ is a progressively measurable $\ell^{2}$ valued process
then this is also the case for $A^{-1/2}\jm{-s}y,$ according to Lemma \ref{measurability 2}.
With $y=\jm{-s}x,$ it follows that $A^{-1/2}x$ is progressively measurable and then
from Lemma \ref{measurability 2} that $\eta$ given by (\ref{prtf eq 3}) is  progressively measurable.
Let $\theta^{1}$ be given by equation (\ref{prtf eq 4}).
Using now  that $\sequiv$ is unitary, that
$\|\theta^{1}_{t}\|_{\dualh{}{}} \leq \|\l_{t}/\zcpxd{t}\|_{\multh{}{}} \|\ \sequiv \eta_{t}\|_{\dualh{}{}}$
and using Theorem \ref{Th price} and  H\"older's inequality, it follows that
\begin{equation} \label{proof th main 4}
\theta^{1} \in \derprod{}{0}(L^{2}(\mathbb{T},\dualh{}{})),
\end{equation}
where $\theta^{1}$ is progressively measurable. Since $\theta^{1}$ satisfies equation (\ref{prtf eq 1})
by construction and formulas (\ref{proof th main 1}) and (\ref{proof th main 4}) hold,
the hypotheses of Lemma \ref{Lm sf prtf} are satisfied, so $\theta^{1} \in \prtfs.$
It also follows that the mapping
$\derprod{}{s} \ni X \mapsto \theta^{1} \in \prtfs \cap \derprod{}{0}(L^{2}(\mathbb{T},\dualh{}{}))$
is continuous.

We define $a$ as in formula (\ref{proof th main 5}). Lemma \ref{Lm sf prtf} then gives
that $\theta^{0} \in \prtfs \cap \derprod{}{0}(L^{2}(\mathbb{T},\dualh{}{})),$ that
$\theta \in \sfprtfs$ and that $\theta^{0}$ has the announced continuity property.
\qed \\ \vspace{-2ex}

\noindent \textbf{Proof of Theorem \ref{th converse main}}
Let the hypotheses of the theorem be satisfied and let
$\theta \in \sfprtfs\cap \derprod{}{0}(L^{2}(\mathbb{T},\dualh{}{})).$
According to Theorem \ref{Th mart meas}, $\prtfpxd{\timeh}(\theta) \in L^{2}(\Omega,Q,\mathcal{F}).$
The self-financing condition (\ref{SFprtf-q}) and Lemma \ref{Lm L2 isomrphism}, show that
$\prtfpxd{\timeh}(\theta) =\uequiv(c,x),$ where $c \in \mathbb{R}$ and
$x \in L_{a}^{2}$ is given by formula (\ref{prtf eq 1}).
Obviously $c \in \derprod{}{s},$ so we only have to prove that $\uequiv(0,x) \in \derprod{}{s}.$
By the construction of $\derprod{}{s}$ it is enough to prove that $x \in L_{a,s}^{p},$ for all $p \geq 2,$ 
which is equivalent to that $x$ is progressively measurable and $\processjm{s}x \in L^{p}(Q,L^{2}(\mathbb{T},\ell^{2} )),$
for all $p \geq 2,$ where $\processjm{}$ is given by formula (\ref{J}). As $x$ is progressively measurable,
Lemma \ref{Lm D0 P Q} shows that
it is  sufficient to check that $\processjm{s}x \in \derprod{}{0}(L^{2}(\mathbb{T},\ell^{2})).$

For the moment let us suppose that a.e.  $(t,\omega),$
\begin{equation}
 \jm{s}(B_{t}(\omega))' \in L(\ell^{2},\dualh{}) \quad \text{and} \quad \|\jm{s}(B_{t}(\omega))'\| \leq k(\omega),
\label{poof convers 1}
\end{equation}
where the norm is the operator norm. Since $(b_{t}(\omega))'=(B_{t}(\omega))'\bar{q}_{t}(\omega),$
where $\bar{q}_{t}(\omega)=\zcpxd{t}(\omega)/l_{t},$ it follows from (\ref{poof convers 1}) that
$\|\jm{s}x_{t}(\omega)\|_{\ell^{2}}=\| \jm{s}(B_{t}(\omega))'\bar{q}_{t}(\omega)\theta_{t}(\omega)\|_{\ell^{2}}$
        $ \leq k(\omega) \|\bar{q}_{t}(\omega)\theta_{t}(\omega)\|_{\dualh{}}.$
Using that $\|Gf\|_{\dualh{}} \leq C \|G\|_{\multh{}}\|f\|_{\dualh{}},$ where
$C$ only depends on $\sobolevorder,$  we obtain
$$\|\jm{s}x_{t}(\omega)\|_{\ell^{2}} \leq k(\omega) \|\bar{q}_{t}(\omega)\|_{\multh{}}
          \|\theta_{t}(\omega)\|_{\dualh{}}.$$
H\"older's inequality, with $1/p=1/p_{1}+1/p_{2}+1/p_{3},$ $p <p_{1},p_{2},p_{3} < \infty,$ gives
\begin{equation}
\|\processjm{s}x\|_{L^{p}(P,L^{2}(\mathbb{T},\ell^{2} ))} \leq C \|k\|_{L^{p_{1}}(P)}
 \|\bar{q}\|_{L^{p_{2}}(P,L^{\infty}(\mathbb{T}, \multh{}))}
     \|\theta\|_{L^{p_{3}}(P,L^{2}(\mathbb{T}, \dualh{}))}.
\label{poof convers 3}
\end{equation}
Since by hypothesis $k \in \derprod{}{0}$  and  $\theta \in \derprod{}{0}(L^{2}(\mathbb{T},\dualh{}{})),$
the norms of  $k$ and $\theta$ on the right hand side of (\ref{poof convers 3}) are finite.
The norm of $\bar{q}$ is also finite, according to Theorem \ref{Th price}, so
$\processjm{s}x \in L^{p}(P,L^{2}(\mathbb{T},\ell^{2} )) ,$ for all $p \geq 2.$
This proves that $\processjm{s}x \in \derprod{}{0}(L^{2}(\mathbb{T},\ell^{2})).$

It remains to prove (\ref{poof convers 1}). If $x \in \ell^{s,2},$ then it follows from
the definition of the process $A$ and the hypothesis of the theorem that
$\|B_{t}(\omega)\jm{s}x\|_{\h{}}^{2}=(\jm{s}x,(B_{t}(\omega))^{*}B_{t}(\omega)\jm{s}x)_{\h{}}
= \|(A_{t}(\omega))^{1/2}\jm{s}x \|_{\ell^{2}}^{2} \leq (k(\omega))^{2} \|x\|_{\ell^{2}}^{2}.$
$B_{t}(\omega)\jm{s}$ from $\ell^{2}$ to $\h{}$ is then closeable and its closure
$K_{t}(\omega) \in L(\ell^{2},\h{})$ has norm bounded by $k(\omega).$ We have
$(K_{t}(\omega))^{*}=\jm{s}(B_{t}(\omega))^{*},$ since $B_{t}(\omega) \in L(\ell^{2},\h{}).$
This shows that $ \jm{s}(B_{t}(\omega))^{*} \in L(\ell^{2},\h{})$ has norm bounded
by $k(\omega).$ The relation $(B_{t}(\omega))'=(B_{t}(\omega))^{*}\sequiv^{-1},$
where $\sequiv$ is the isomorphism defined in (\ref{S equiv}), now gives (\ref{poof convers 1}).
\qed \\ \vspace{-2ex}

\noindent \textbf{Proof of Theorem \ref{th opt port}}
We only consider the case of $u'>0,$ since the case of $u'(x)=0$ for some $x$ is so similar.
Let the hypotheses of the theorem be satisfied. According to Corollary 3.4 of
reference \cite{I.E.-E.T}, the portfolio $\hat{\theta}$ is a solution of
equation (\ref{opt prob}), if $\prtfpxd{\timeh}(\hat{\theta})=\hat{X},$
where  $\hat{X}=\varphi (\lambda \xi_{\timeh})$ for a certain $\lambda >0.$
$\varphi $ is $C^{1}$ and $\ln(\xi_{\timeh}) \in \derprod{1}{s},$ so
$\jm{s}\mder{t}\hat{X}=\lambda \xi_{\timeh}\varphi' (\lambda \xi_{\timeh})\jm{s}\mder{t}\ln(\xi_{\timeh}).$
This gives
$\|\processjm{s}\mder{}\hat{X}\|_{L^{2}( \mathbb{T},\ell^{2})}
  =|\lambda \xi_{\timeh}\varphi' (\lambda \xi_{\timeh})| \, \|\processjm{s}\mder{}\ln(\xi_{\timeh})\|_{L^{2}( \mathbb{T},\ell^{2})}.$
Inequality (\ref{utility cond 2}) gives
$\|\processjm{s}\mder{}\hat{X}\|_{L^{2}( \mathbb{T},\ell^{2})}
\leq C ((\lambda \xi_{\timeh})^{p}+(\lambda \xi_{\timeh})^{-p}) \|\processjm{s}\mder{}\ln(\xi_{\timeh})\|_{L^{2}( \mathbb{T},\ell^{2})}.$
Theorem \ref{Th mart meas} shows that $(\lambda \xi_{\timeh})^{p}+(\lambda \xi_{\timeh})^{-p} \in L^{q}(\Omega, P),$
for all $q \geq 1.$ By hypothesis
$\|\processjm{s}\mder{}\ln(\xi_{\timeh})\|_{L^{2}( \mathbb{T},\ell^{2})} \in \derprod{}{0},$
so H\"older's inequality now gives that
$\|\processjm{s}\mder{}\hat{X}\|_{L^{2}( \mathbb{T},\ell^{2})} \in \derprod{}{0},$
i.e. $\mder{}\hat{X} \in \derprod{}{0}(L^{2}( \mathbb{T},\ell^{s,2})).$
By Theorem 3.3 of reference \cite{I.E.-E.T}, $\hat{X} \in \derprod{}{0}.$
Corollary \ref{cor inclusion} then gives that $\hat{X} \in \derprod{1}{s}.$
We can now apply Theorem \ref{th main}, which proves the existence of $\hat{\theta}.$
\qed
\appendix
\section{Auxiliary results}
\label{App}
\noindent \textbf{Proof of Lemma \ref{Lm L2 isomrphism}}
We first prove that the mapping $(c,x) \mapsto X=\uequiv(c,x)$  is  unitary on
$\mathbb{R} \oplus L_{a}^{2}$ to $L^{2}(\Omega ,Q,\mathcal{F}).$ %

The operator $\uequiv$ is isometric, so its range is a closed subspace of $L^{2}(\Omega ,Q,\mathcal{F}).$
In fact (cf. Proposition 4.13 of \cite{DaPrato-Zabczyk}),
$\|\uequiv(c,x)\|_{L^{2}(Q)}^{2}=c^{2}+\|x\|_{L_{a}^{2}}^{2}.$
It is sufficient to prove that $\uequiv$ has dense range. Let $h \in L^{2}( \mathbb{T},\ell^{2})$ and let
$$\mathcal{E}_{t}(h)=\exp{ \left( -\frac{1}{2}\int_{0}^{t}\sum_{i \in \mathbb{N}^{*}} (h(s,i))^{2}ds
    +\int_{0}^{t}\sum_{i \in \mathbb{N}^{*}}h(s,i)  d\wienerq{i}{s})\right)}.$$
$\mathcal{E}_{\timeh}(h)$ is in the range of $\uequiv,$ since $\mathcal{E}(h)h \in L_{a}^{2}$
and by It\^o's lemma (Theorem 4.17 of \cite{DaPrato-Zabczyk}):
$$\mathcal{E}_{\timeh}(h)=1+\int_{0}^{\timeh}\sum_{i \in \mathbb{N}^{*}}\mathcal{E}_{s}(h)h(s,i)  d\wienerq{i}{s}.$$
We have $L^{2}( \mathbb{T},\ell^{2}) \cong L^{2}( \mathbb{T} \times \mathbb{N}^{*})$
and the measure $\nu$ is atomless. The linear span of
$\{\mathcal{E}_{\timeh}(h) \, | \, h \in L^{2}( \mathbb{T},\ell^{2}) \}$
is then dense in $L^{2}(\Omega ,Q,\mathcal{F})$ (cf. Lemma 1.1.2 of \cite{Nual71}),
which proves that $\uequiv$ is a unitary operator.

To prove the second part of the lemma we fix $p \geq 2.$
For $(c,x) \in \mathbb{R} \oplus L_{a}^{2}$ let $X=\uequiv(c,x),$ for $0 \leq t \leq \timeh$ let
\begin{equation} \notag %
Y_{t}=\sum_{i \in \mathbb{N}^{*}}\int_{0}^{t}x_{s}^{i}d\wienerq{i}{t}
\end{equation}
and let $Z=\sup_{0 \leq t \leq \timeh }|Y_{t}|.$ In the sequel of this proof $C, C_{1},C_{2}, \ldots $
are positive constants independent of $X$ and $(c,x).$
 Applying the BDG inequalities we obtain
$$\|X\|_{L^{p}(Q)} \leq |c|+\|Z\|_{L^{p}(Q)} \leq |c|+ C\|x\|_{L_{a}^{p}}.$$
This shows that
\begin{equation}  \label{proof isomrphism 1}
\uequiv(\mathbb{R} \oplus L_{a}^{p}) \subset L^{p}(\Omega ,Q,\mathcal{F}).
\end{equation}
Given $X \in L^{p}(\Omega ,Q,\mathcal{F}),$ then $X \in L^{2}.$  By the first part of
the lemma it follows that $(c,x) =\uequivadj X \in \mathbb{R} \oplus L_{a}^{2}.$ Since $\uequivadj$ is continuous,
$|c| \leq C_{1}\|X\|_{L^{2}(Q)} \leq C_{1} \|X\|_{L^{p}(Q)}.$ The BDG inequalities give
$\|x\|_{L_{a}^{p}} \leq C_{2} \|Z\|_{L^{p}(Q)}.$
Applying  Doob's $L^{p}$ inequalities and using that $|Y|^{p}$ is a submartingale, we obtain that
$$\|x\|_{L_{a}^{p}}  \leq C_{3} \sup_{0 \leq t \leq \timeh }\|Y_{t}\|_{L^{p}(Q)}  \leq C_{3} \|X\|_{L^{p}(Q)}.$$
This proves that $\uequivadj L^{p}(\Omega ,Q,\mathcal{F}) \subset \mathbb{R} \oplus L_{a}^{p}.$
Since $\uequiv$ is unitary it follows that
$L^{p}(\Omega ,Q,\mathcal{F}) \subset \uequiv (\mathbb{R} \oplus L_{a}^{p}),$
which together with (\ref{proof isomrphism 1}) proves that
$\uequiv(\mathbb{R} \oplus L_{a}^{p}) = L^{p}(\Omega ,Q,\mathcal{F}).$
This gives by unitarity
$\uequivadj L^{p}(\Omega ,Q,\mathcal{F}) = \mathbb{R} \oplus L_{a}^{p}.$

Finally the restriction $B \in L(\mathbb{R} \oplus L_{a}^{p},L^{p}(\Omega ,Q,\mathcal{F}))$
of $\uequiv$ is a homeomorphism since $B^{-1}$ is the restriction of $\uequivadj$
to $L^{p}(\Omega ,Q,\mathcal{F}).$
\qed \\ \vspace{-2ex}

In the sequel $E,$ $E_{1}$ and $E_{2}$ are separable Hilbert spaces.
The next lemma collects some well-known results on polar decomposition, cf. Ch VI, \S7 of \cite{Kato66}.
We recall that, if %
$K$ is a densely defined closed operator from $E_{1}$  to $E_{2}$ with
adjoint $K^{*},$ then according to von Neumann's theorem, $K^{*}K$ is a positive self-adjoint operator
in $E_{1}.$ Its positive square-root is then well-defined.
\begin{lemma} \label{range}
Let $E_{1}$ and $E_{2}$ be Hilbert spaces and let $K$ be a densely defined closed
operator from $E_{1}$  to $E_{2}.$ The following statements are true:
$i)$ $\mathcal{R}(K^{*})=\mathcal{R}((K^{*}K)^{1/2})$ and $\mathcal{K}(K)=\mathcal{K}((K^{*}K)^{1/2}),$
$ii)$ If $\mathcal{K}(K)=\{0\},$ then $\mathcal{D}((K^{*}K)^{-1})$ is dense in $E_{1},$
$\mathcal{D}((K^{*}K)^{-1}) \subset \mathcal{D}((K^{*}K)^{-1/2})$
and the closure  of $K(K^{*}K)^{-1/2}$ is an isometric operator $S \in L(E_{1},E_{2}),$
$iii)$ If $\mathcal{K}(K)=\{0\}$ and $x \in \mathcal{D}((K^{*}K)^{-1/2}),$ then
$K^{*}S(K^{*}K)^{-1/2}x=x.$
\end{lemma}
\textbf{Proof:} Let $D=K^{*}K.$

$i)$ This statement follows from Problem 2.33, \S7, Ch. VI of \cite{Kato66}.

$ii)$ $\mathcal{K}(D)=\mathcal{K}(K)=\{0\}.$ Since $D$ is selfadjoint it
follows that $\mathcal{D}(D^{-1})$ is dense in $E_{1}.$ Using the spectral
resolution of $D$ (cf. Ch XI \S12 \cite{Yosida}), we obtain
$\mathcal{D}(D^{-1}) \subset \mathcal{D}(D^{-1/2}).$ Let $x \in \mathcal{D}(D^{-1/2}) \cap \mathcal{D}(D^{1/2}).$
Then $\|KD^{-1/2}x\|_{E_{2}}^{2}=(D^{-1/2}x,K^{*}KD^{-1/2}x)_{E_{2}}=\|x\|_{E_{1}}^{2}.$
Since $\mathcal{D}(D^{-1/2}) \cap \mathcal{D}(D^{1/2})$ is dense in $E_{1},$ it now follows that the closure $S$
is an isometric operator.

$iii)$ Let $x \in \mathcal{D}(D^{-1}).$ Then $D^{-1/2}x \in \mathcal{D}(D^{-1/2})$
and $S=KD^{-1/2}$ on $\mathcal{D}(D^{-1/2}),$ so  $K^{*}SD^{-1/2}x=K^{*}KD^{-1}x=x.$
This equality extends by continuity to $x \in \mathcal{D}(D^{-1/2}).$
\qed \\ \vspace{-2ex}

\noindent The spectrum $\sigma(K),$ of a compact selfadjoint $K$ operator on $E,$ is
real, denumerable and zero is the only possible accumulation point.
The spectral resolution of $K$ is given by
\begin{equation}
K=\sum_{\lambda \in \sigma(K)}\lambda \, e_{K}( \{\lambda \}),
\label{spect rep 1}
\end{equation}
where $e_{K}$ is the corresponding resolution of the identity defined on the Borel
subsets of $\mathbb{R}.$ If $f$ is a real valued function on $\mathbb{R},$ then the operator $f(K)$ in $E$
is given by
\begin{equation}
f(K)=\sum_{\lambda \in \sigma(K)} f(\lambda) \, e_{K}( \{\lambda \}),
\label{spect rep 2}
\end{equation}
on its domain $\mathcal{D}(f(K))= \{x \in E \; | \; \sum_{\lambda \in \sigma(K)} |
        f(\lambda)|^{2} \, \|e_{K}( \{\lambda \})x\|_{E}^{2} < \infty \}.$
\begin{lemma} \label{measurability 1}
Let $A$ be the set of compact selfadjoint operators in $E,$ endowed with the operator norm topology.
If $f: \mathbb{R} \rightarrow \mathbb{R}$ is a Borel function, %
bounded on bounded subsets of $\mathbb{R},$
then the function $A \times E \ni (K,x) \mapsto f(K)x \in E$ is Borel measurable.
Moreover the mapping $\mathbb{R} \times A \times E \ni (\lambda, K,x) \mapsto e_{K}( \{\lambda \})x \in E$
is Borel measurable.
\end{lemma}
\textbf{Proof:} Here $I \in L(E)$ is the identity operator and $L(E)$ is given the operator norm topology.
Let $B=\{K- \lambda I \, | \, K \in A, \, \lambda \in \mathbb{R}\}$ be endowed with the operator norm topology.
$B$ is a closed subalgebra of $L(E).$
The subspace  $A_{0}=\{K- \lambda I \in B \, | \, \lambda \neq 0\}$ is open in $B.$
For given $M=K- \lambda I \in A_{0},$ the space $\mathcal{K}(M)$ has finite dimension, since $K$ is compact,
and $\mathcal{R}(M)$ is a closed subspace of $E.$ It now follows as in the finite dimensional case
(cf. Chap. $1,$ Lemma $4.4$ of \cite{K-S 99}), that the mapping $A_{0} \times E \ni (M,x) \mapsto e_{M}( \{0 \})x \in E$
is Borel measurable. Since $(\mathbb{R} -\{0\}) \times A \ni (\lambda,K) \mapsto K- \lambda I \in A_{0}$
is continuous and $e_{K-\lambda I}( \{0 \})=e_{K}( \{\lambda \}),$ the mapping
$F: (\mathbb{R} -\{0\}) \times A \times E \rightarrow E,$ where $F( \lambda, K, x) = e_{K}( \{\lambda \})x,$
is Borel measurable.

Suppose that $f$ satisfies the hypothesis of the lemma and let $G(K,x)=f(K)x.$
We first consider the case of $f(x) =0$ for all $x \leq 0.$
For $K \in A,$ let $\mu_{1}(K) \geq \cdots \mu_{n}(K)  \cdots \geq 0$ be the decreasing sequence
of positive eigenvalues of $K,$ each repeated a number of times equal to the multiplicity
of the  eigenvalue. The function $A \ni K \mapsto \mu_{n}(K)$ is then continuous (cf. \cite{Kato66}, Ch. IV, \S3.5).
Define $\mu_{0}(K)=\mu_{1}(K)+1$ and $u: \mathbb{R}^{2} \rightarrow \mathbb{R}$
by $u(x,y)=0$ if $x \leq y$ and  $u(x,y)=1$ if $x > y.$ $u$ is a Borel function.
It follows from (\ref{spect rep 2}) that
\begin{equation}
G(K,x)=\sum_{n=1}^{\infty} u(\mu_{n-1}(K),\mu_{n}(K)) f(\mu_{n}(K)) \, e_{K}( \{\mu_{n}(K) \})x, \; x \in E.
\label{proof measurability 1}
\end{equation}
The mapping $G: A \times E \rightarrow E$ is Borel measurable. In fact, by the continuity of $\mu_{n}$
and the measurability of $u$ and $F,$ each term in the sum (\ref{proof measurability 1}) is Borel measurable,
as a function of $(K,x).$ The sum (\ref{proof measurability 1}) converges pointwise $(K,x)$ in $E$
to $G(K,x),$ so $G$ is Borel measurable (cf. Theorem 5.6.3 of \cite{SCHW 93}).

Next we consider the case of $f(x) =0$ for all $x \geq 0.$ Similarly as to the previous case,
it follows that $G$ is a Borel function. Finally we consider the case of $f(x)=0$ for $x \neq 0$ and $f(0)=a.$
From the two previous cases it follows that $(K,x) \mapsto h(K)x$ is Borel measurable, where
$h(0)=0$ and $h(x)=1$ for $x \neq 0.$ Since $G(K,x)=x-ah(K)x,$ it follows that $G$ is a Borel function.
The case of a general $f$ now follows by the decomposition $f=f_{-}+f_{0}+f_{+},$
where the support of $f_{-},$ $f_{0}$ and $f_{+}$ is a subset of $]-\infty,0[,$ $\{0\}$ and $]0, \infty[$ respectively.

To prove the last statement, we note that $A \times E \ni ( K,x) \mapsto e_{K}( \{0 \})x \in E$ 
is Borel measurable, which follows from the identity
$e_{K}( \{0 \})x=x-h(K)x,$ where $h(0)=1,$ $h(\lambda)=0$ for $\lambda \neq 0.$
The measurability of $\mathbb{R} \times A \times E \ni (\lambda, K,x) \mapsto e_{K}( \{\lambda \})x \in E$
now follows from the Borel measurability of $F.$
\qed
\begin{lemma} \label{measurability 2}
Let $A$ be the set of compact  operators  with trivial kernel from $E_{1}$ to $E_{2},$
endowed with the operator norm topology.
If $S_{K}$ is the closure of the operator $K(K^{*}K)^{-1/2}$
then the function $A \times E_{1} \ni (K,x) \mapsto S_{K}x \in E_{2}$ is Borel measurable.
Moreover, if $L \in L(E_{1})$ and $A'$ is the subspace of elements  $K \in A$ such that
$\mathcal{R}(L) \subset \mathcal{D}((K^{*}K)^{-1/2}),$ then $A' \times E_{1} \ni (K,x) \mapsto (K^{*}K)^{-1/2}Lx \in E_{1}$ is Borel measurable.
\end{lemma}
\textbf{Proof:} Let $f_{n}(x)=\sqrt x$ if $x \geq 1/n$ and $f_{n}(x)=0$ if $x< 1/n,$ for $n \in \mathbb{N}^{*}.$
The function $A \ni K \mapsto K^{*}K \in A$ is continuous and $K^{*}K$ is selfadjoint.
Let $F_{n}(K,x)=Kf_{n}( K^{*}K)x.$ Lemma \ref{measurability 1} shows that $F_{n}: A \times E_{1} \rightarrow E_{2}$
is Borel measurable. Since $F_{n}(K,x)$ converges pointwise to $S_{K}x$ in $E_{2}$
as $n \rightarrow \infty,$ it follows that $(K,x) \mapsto S_{K}x$ is Borel measurable.
To prove the second statement, we note that $A' \times E_{1} \ni (K,x) \mapsto f_{n}( K^{*}K)Lx \in E_{1}$
is Borel measurable. Since $Lx \in \mathcal{D}((K^{*}K)^{-1/2}),$ the sequence $f_{n}( K^{*}K)Lx$
converges pointwise $(K,x)$ in $E_{1}$ to $(K^{*}K)^{-1/2}Lx.$ It follows that $(K,x) \mapsto (K^{*}K)^{-1/2}Lx $
is Borel measurable.
\qed

We shall define two mappings, $g_{0}$ and  $g_{1},$ on the space of selfadjoint Hilbert-Schmidt operators on $E.$ 
They will satisfy $g_{0}(K) \in \mathcal{K}(K)$ and $g_{1}(K) \in \mathcal{R}(K)^{c}.$
Let $\{u_{n}\}_{n \geq 1} $ be an orthonormal basis in $E$ and let $K$ be a selfadjoint Hilbert-Schmidt operators on $E.$
We define
\begin{equation}
g_{0}(K)=0 \;  \text{if}  \; \mathcal{K}(K)=\{0\}  \; \text{and}  \; g_{0}(K)=\frac{e_{K}( \{0 \})u_{N(K)}}{\|e_{K}( \{0 \})u_{N(K)}\|} \;
 \text{if}  \; \mathcal{K}(K) \neq \{0\},
\label{g0}
\end{equation}
where $N(K)= \min \{n \, | \, e_{K}( \{0 \})u_{n} \neq 0 \}.$
If $\lambda \notin \sigma (K),$  then let $h(K,\lambda)=0$ and if $\lambda \in \sigma (K)$ has multiplicity $m,$ then let
$$h(K,\lambda)=v_{1}+\cdots+v_{m},$$
where $\{v_{1},\cdots,v_{m}\}$ is the orthonormal basis in $e_{K}( \{\lambda \})E$ given by the
Schmidt orthonormalization of $\{e_{K}( \{\lambda \})u_{n}\}_{n \geq 1}.$
More precisely let $P_{0}=e_{K}( \{\lambda \})$ and we construct
inductively $n_{1},\ldots,n_{m},$ $v_{1},\ldots,v_{m}$ and $P_{1},\ldots,P_{m}$ by:

$n_{1}=\min \{n \, | \, P_{0}u_{n} \neq 0 \},$
$v_{1}=P_{0}u_{n_{1}}/\|P_{0}u_{n_{1}}\|$ and $P_{1}x=(v_{1},x)_{E}v_{1}.$
$Q_{k+1}=P_{0}-\sum_{i=1}^{k}P_{i},$
$n_{k+1}=\min \{n \, | \, Q_{k+1}u_{n} \neq 0 \},$ $v_{k+1}=Q_{k+1}u_{n_{k+1}}/\|Q_{k+1}u_{n_{k+1}}\|$
and $P_{k+1}x=(v_{k+1},x)_{E}v_{k+1}.$

We now define $g_{1}$ by
\begin{equation}
g_{1}(K)=g'_{1}(K)/\|g'_{1}(K)\|, \; \text{where} \; g'_{1}(K)=\sum_{\lambda \in \sigma (K)}\lambda h(K,\lambda)+g_{0}(K).
\label{g1}
\end{equation}
\begin{lemma} \label{g_{0} g_{1}}
Let $A$ be the set of selfadjoint Hilbert-Schmidt operators on $E,$
endowed with the operator norm topology. The maps $g_{i}: A \rightarrow E,$ $i=0,1$ given by (\ref{g0}) and (\ref{g1})
are Borel measurable.
For every $K \in A$ %
the following two properties are satisfied:
$i)$ $g_{0}(K) \in \mathcal{K}(K)$ and if $\mathcal{K}(K) \neq 0$ then  $\|g_{0}(K)\|=1.$
$ii)$ $g_{1}(K) \notin \mathcal{R}(K)$ and $\|g_{1}(K)\| =1.$
\end{lemma}
\textbf{Proof:} Since Hilbert-Schmidt operators are compact it follows from Lemma \ref{measurability 1} that
$\mathbb{R} \times A \times E \ni (\lambda, K,x) \mapsto e_{K}( \{\lambda \})x \in E$
is Borel measurable. The measurability of $g_{0}$ then follows from that $N$ is measurable and that $x \mapsto x/\|x\|$
is measurable on $E-\{0\}.$ Similarly, for given $(K,\lambda),$ $v_{i}$ is a measurable function of a finite number
of the variables $e_{K}( \{\lambda \})u_{n}.$ Therefore $h:  A \times \mathbb{R} \rightarrow E$ is measurable.
The sum in (\ref{g1}) converges. In fact, using that $\|h(K,\lambda)\|_{E}^{2}$ is equal to the dimension
of $e_{K}( \{\lambda \})E,$ it follows that
\begin{equation}
\|g'_{1}(K)-g'_{0}(K)\|_{E}^{2}=\sum_{\lambda \in \sigma (K)}\lambda^{2}\|h(K,\lambda)\|_{E}^{2}=\|K\|_{H-S}^{2}.
\label{g1'}
\end{equation}
The Borel measurability of $g_{1}$ follows from the pointwise convergence.

Statement $i)$ is obvious and we prove $ii)$.
If $\mathcal{K}(K) \neq \{0\},$
then $e_{K}( \{0 \})g_{1}(K)$ $=g_{0}(K) \neq 0$ and $e_{K}( \{0 \})\mathcal{R}(K)=\{0\},$
show that $g'_{1}(K) \notin \mathcal{R}(K).$
Let $\mathcal{K}(K)= \{0\},$ suppose that $g'_{1}(K) \in \mathcal{R}(K),$ let $x \in E$
be the unique element such that $g'_{1}(K)=Kx$ and let $x_{\lambda}=e_{K}( \{\lambda \})x.$
Then $Kx_{\lambda}=e_{K}( \{\lambda \})Kx=Kh(K,\lambda),$ so $x_{\lambda}=h(K,\lambda).$
This gives that $\|x\|_{E}^{2}=\sum_{\lambda \in \sigma (K)}$ $\|h(K,\lambda)\|_{E}^{2}=\infty.$
This is a contradiction, so $g'_{1}(K) \notin \mathcal{R}(K).$
Hence  $g'_{1}(K) \notin \mathcal{R}(K)$ for every $K \in A.$
In particular  $g'_{1}(K) \neq 0,$ for every $K \in A,$ so $g_{1}(K)$ is well-defined
and $\|g_{1}(K)\| =1.$
\qed
\bibliographystyle{amsplain}

\begin{thebibliography}{99}
%
%
%
\bibitem{Bj-Ka-Ru97}
Bj\"ork, T., Kabanov, Y. and Runggaldier, W.: {\em Bond market structure in the presence
of marked point processes}, Mathematical Finance, {\bf 7}, 211--239 (1997).
\bibitem{Bj-Ma-Ka-Ru97}
Bj\"ork, T., Masi, G., Kabanov, Y. and Runggaldier, W.: {\em Toward a general theory
of bond markets}, Finance and Stochastics, {\bf 1}, 141--174 (1997).
%
%
%
\bibitem{Calderon}
Calderon, A.P.: {\em Lebesgue spaces of differentiable functions and distributions},
Proc. Symp. Pure Math. IV, AMS 1961, 33--49.
\bibitem{Carmona-Tehr}
Carmona, R. and Tehranchi, M.: {\em A Characterization of Hedging Portfolios for
Interest Rate Contingent Claims}, Preprint March 24, 2003.
\bibitem{DaPrato-Zabczyk}
Da Prato, G. and Zabczyk, J.: {\em Stochastic Equations in Infinite Dimensions},
Encyclopedia of Mathematics and its Applications, Cambridge University Press, 1992.
\bibitem{DeDonno Pratelli}
De Donno, M. and Pratelli, M.: {\em On the use of measure-valued strategies in bond markets},
Finance and Stochastics, {\bf 8}, 87--109 (2004).
\bibitem{I.E.-E.T} Ekeland, I. and Taflin, E.: {\em A Theory of Bond Portfolios},
%
To appear in Annals of Applied Probability; \\ \texttt{http://arxiv.org/abs/math.OC/0301278}
%
%
%
%
%
%
%
%
\bibitem{HJM92} Heath, D.C., Jarrow, R.A. and Morton, A.: {\em Bond pricing and the
term structure of interest rates: a new methodology for contingent claim valuation},
Econometrica, {\bf 60}, 77--105 (1992).
\bibitem{Horm}  H\"ormander, L.: {\em The analysis of  linear  partial  differential
operators}, Vol. I, Springer-Verlag 1985.
\bibitem{K-S 99}
Karatzas, I. and Shreve, S.E.: {\em Methods of Mathematical Finance}, Applications of Mathematics, Volume 9,
Springer-Verlag 1999.
\bibitem{Kato66}
Kato, T. {\em Perturbation Theory for Linear Operators}, Die Grundleheren der
mathematischen Wissenschaften, Volume 132, Springer-Verlag, New York 1966.
%
%
%
%
%
%
%
%
%
%
%
%
%
%
%
%
%
%
%
\bibitem{Nual71}
Nualart D.: {\em The Malliavin Calculus and Related Topics}, Probability and its Applications,
Springer-Verlag, 1991.
\bibitem{Pham 2003}
Pham, H.: {\em A predictable decomposition in infinite asset model with jumps.
Application to hedging and optimal investment}, Stochastics and Stochastic Reports, {\bf 5}, 343--368 (2003).
%
%
%
\bibitem{Revuz-Yor}
Revuz, D. and Yor, M.: {\em Continuous Martingales and Brownian Motion},
Grundlehren der mathematischen Wissenschaften,
Band 293, Springer-Verlag.
%
%
%
%
%
%
%
\bibitem{SCHW 93}
Schwartz, L.: {\em Analyse III}, Hermann 1993, Paris.
\bibitem{Yosida}
Yosida, K.: {\em Functional Analysis}, Grundlehren der mathematischen Wissenschaften,
Band 123, Springer-Verlag.
\end{thebibliography}

\end{document}